\documentclass[11pt]{article}

\usepackage[margin=1in]{geometry}
\usepackage[T1]{fontenc}
\usepackage{lmodern}
\usepackage{amsmath,amssymb,amsthm,amsfonts,mathtools}
\usepackage{booktabs,array}
\usepackage[numbers,square,sort&compress]{natbib}
\usepackage{microtype}
\usepackage{needspace}
\usepackage[colorlinks=true,linkcolor=blue,citecolor=blue,urlcolor=blue]{hyperref}

\numberwithin{equation}{section}
\allowdisplaybreaks
\theoremstyle{plain}
\newtheorem{theorem}{Theorem}[section]
\newtheorem{proposition}[theorem]{Proposition}
\newtheorem{corollary}[theorem]{Corollary}
\newtheorem{lemma}[theorem]{Lemma}
\theoremstyle{definition}
\newtheorem{definition}[theorem]{Definition}
\newtheorem{example}[theorem]{Example}
\theoremstyle{remark}
\newtheorem{remark}[theorem]{Remark}

\newcommand{\Lr}{\operatorname{L}}
\newcommand{\OK}{\mathcal O_K}
\newcommand{\Q}{\mathbb Q}
\newcommand{\Z}{\mathbb Z}
\newcommand{\R}{\mathbb R}

\newcommand{\Bx}{\mathcal B}
\newcommand{\vp}{v_{\mathfrak p}}

\title{\textbf{Two-term Rogers dilogarithm identities over cubic fields}}
\author{Cetin Hakimoglu-Brown}
\date{November 22, 2026}
\hypersetup{
  pdftitle={Two-term Rogers dilogarithm identities over cubic fields},
  pdfauthor={Cetin Hakimoglu-Brown},
  pdfsubject={Rogers dilogarithm identities and exceptional units in cubic fields},
  pdfkeywords={Rogers dilogarithm, exceptional units, cubic fields, unit equations}
}

\begin{document}
\maketitle

\begin{abstract}
We study normalized two-term Rogers dilogarithm relations in cubic fields. For integral,
nonunimodular rank-two monomial configurations whose base argument generates
a totally real cubic field, finite-prime valuations force the arguments
and their complements to be units. Normalization preserves the generated
order and reduces the problem to two cubic families without a maximal-order
assumption. In the simplest cubic family, an existing unit-equation
classification yields exactly ten relations up to simultaneous reflection,
with distinct coefficients 2, 3, 7, 13, 19, 31 and 97. Explicit five-term
certificates prove all ten evaluations. In the second family, a sharpened lower bound for the
height of a solution, together with the two-logarithm theorem of Laurent, Mignotte and
Nesterenko, reduces the unit equation to a short certified computation, and the family
contributes no relation in this class. A unit-index argument restricts, independently, wedge relations between
arbitrary exceptional units of non-Galois totally real cubic fields: any such relation would have a coefficient with no
prime factor below 23 and would contradict Ennola's conjecture on fundamental units.
In cubic fields with one real embedding, the only five-term certified relation
between exceptional units with coefficient greater than one is Lewin's identity.
Together, these results give exactly eleven identities in the two classes
considered, up to the stated reflections.
\end{abstract}



\noindent\begin{minipage}{\textwidth}
\small\tableofcontents
\end{minipage}
\vspace{1em}

\section{Introduction}\label{sec:intro}

For $0<x<1$, the Rogers dilogarithm is
\begin{equation}\label{eq:Ldef}
\Lr(x)=\operatorname{Li}_2(x)
+\frac12\log x\,\log(1-x).
\end{equation}
It satisfies the five-term identity
\begin{equation}\label{eq:fiveterm}
\Lr(x)+\Lr(y)
=
\Lr(xy)
+\Lr\!\left(\frac{x(1-y)}{1-xy}\right)
+\Lr\!\left(\frac{y(1-x)}{1-xy}\right),
\qquad 0<x,y<1,
\end{equation}
and the reflection identity
\begin{equation}\label{eq:reflection}
\Lr(x)+\Lr(1-x)=\Lr(1)=\frac{\pi^2}{6}.
\end{equation}
We consider two-term evaluations
\begin{equation}\label{eq:twoterm}
m\,\Lr(z_0)+\Lr(z_1)=c\pi^2,
\qquad
z_0,z_1\in(0,1)\cap\overline{\Q},
\quad m\in\Z_{>1},\quad c\in\Q.
\end{equation}

Classical examples include Watson's conductor-$7$ identities
\citep{Watson}, Lewin's plastic identity
\citep{Lewin82}, and the two-term identities obtained from
thermodynamic Bethe ansatz equations by \citet{Bytsko}.  A modulus-$9$
calculation of \citet{Kanade} led to the identity appearing as row~3
of Table~\ref{tab:identities}.

\begin{definition}[Five-term certified relation]\label{def:certified}
A Rogers-dilogarithm relation is \emph{five-term certified} if, after exact
algebraic identifications among its arguments, it is obtained from a finite
rational linear combination of the five-term identity
\eqref{eq:fiveterm} and the reflection identity \eqref{eq:reflection}.
\end{definition}

The purpose of the present paper is not to classify every conceivable
numerical equality \eqref{eq:twoterm}.  In the non-degenerate totally real
branch we classify primitive normalized $m:1$ relations supplied by an
integral rank-two monomial configuration and the associated wedge
cancellation; every surviving relation is then proved by an explicit
five-term certificate.  In the rank-one signature-$(1,1)$ branch we classify
the five-term certified relations between exceptional units directly.  A unit
$\varepsilon$ such that $1-\varepsilon$ is also a unit is called an
\emph{exceptional unit} \citep{Nagell64}.

The main arithmetic inputs are a finite-prime rigidity theorem, which forces
the coordinates of a totally real cubic monomial configuration to be units,
and a uniform exclusion of nontrivial exceptional units in the parametrized
orders of the second cubic family. Generated-order normalization connects
these results to the unit-equation classification in the simplest cubic
family. Explicit five-term certificates then establish every surviving
dilogarithm evaluation.

The cubic exceptional-unit classification used throughout the paper is
most conveniently stated in the normalizations used later.

\begin{theorem}[Cubic exceptional-unit classification]
\label{thm:nagell}
Let $K$ be a cubic number field containing an exceptional unit.
\begin{enumerate}
\item If $K$ is totally real, then, up to isomorphism, $K$ is generated
by a root of one of the two normalized families
\[
f_t(X)=X^3-(t+3)X^2+tX+1,
\qquad t\ge-1,
\]
or
\[
g_s(X)=X^3-(s+1)X^2+sX+1,
\qquad s\ge4.
\]
The omitted value $s=3$ gives the conductor-$7$ field already occurring
in the first family.
\item If $K$ has signature $(1,1)$, then
\[
\operatorname{disc}(K)\in\{-23,-31\},
\]
and both fields occur.
\end{enumerate}
Part~(1) is the totally real field-level classification of
\citet{Nagell69}.  Section~\ref{sec:endpoint} recovers the same two
families directly from the elimination polynomial in the positive
normalized presentation.  Part~(2) is the signature-$(1,1)$ part of
\citet{Nagell69}.  Lemma~\ref{lem:generatedorder} strengthens the totally real statement:
an anharmonic transform of the particular unit can be used, preserving its order.
\end{theorem}

The main result can now be stated at the outset.

\begin{theorem}[Main classification theorem]
\label{thm:main}
The following statements hold.
\begin{enumerate}
\item Every non-degenerate normalized integral rank-two configuration
with $K=\Q(x)$ a totally real cubic field is supported on algebraic units.
Consequently its field belongs to one of the two totally real families
in Theorem~\ref{thm:nagell}.
\item In the parametrized simplest-cubic orders, the non-degenerate
normalized $m:1$ relations with $m>1$ arising from wedge cancellation
are exactly the ten identities in Table~\ref{tab:identities}.  Their
distinct coefficients are
\[
2,3,7,13,19,31,97.
\]
\item In the second totally real family, the parametrized orders
$\Z[\theta_s]$, $s\ge4$, have no exceptional units outside the trivial
anharmonic orbit of $\theta_s$, so the family contributes no normalized
$m:1$ relation with $m>1$.  Independently, every relation $m\,\delta(z_0)+\delta(z_1)=0$ with $m\ge2$ between
exceptional units of a non-Galois totally real cubic field has every
prime divisor of $m$ at least $23$, and would contradict Ennola's
conjecture $j_l=1$.
\item In cubic fields of signature $(1,1)$, the five-term certified
relations between exceptional units in $(0,1)$ have no coefficient
$m>1$ in discriminant $-31$, while in discriminant $-23$ the unique
$m>1$ relation, up to simultaneous reflection, is Lewin's identity
\[
2\Lr(\eta)+\Lr(\eta^2)=\frac{\pi^2}{3}.
\]
\end{enumerate}
For the totally real configurations in part~(1), the generated-order
normalization of Lemma~\ref{lem:generatedorder} makes the order-level results
in parts~(2) and (3) sufficient. Consequently, the normalized $m:1$
relations in the totally real integral monomial class consist of the
identities in Table~\ref{tab:identities}.
Together with the exceptional-unit class of part~(4), this gives those ten
relations and Lewin's identity. No assertion here forces arbitrary
signature-$(1,1)$ arguments to be exceptional units.
\end{theorem}

\begin{proof}
Part~(1) is Theorem~\ref{thm:forcing} combined with Lemma~\ref{lem:generatedorder}.  Part~(2) is
Theorem~\ref{thm:identities}; part~(3) is
Theorem~\ref{thm:secondempty} with Corollary~\ref{cor:second}, and
Theorem~\ref{thm:index}; and part~(4) is
Theorem~\ref{thm:lewin}.  The cited results are proved in the sections
below, with the two longest technical arguments and all ten
five-term certificates placed in the appendices.
\end{proof}

The first implication is independent of the dilogarithm.  It states
that a totally real cubic point of
\[
C\cap\phi_M^{-1}(C),\qquad C:\ X+Y=1,
\]
with a non-unimodular integral monomial map $\phi_M$, must be integral
on $C$.  The complete proof is given in
Appendix~\ref{app:forcing}.  The distinction between a five-term certified Rogers relation and an arbitrary
numerical equality is maintained throughout; no converse from numerical
Rogers identities to wedge cancellation or to five-term certification is
assumed.

\section{Bloch groups, the determinant law, and normalized configurations}
\label{sec:setup}

Let $\mathcal P(F)$ be the pre-Bloch group of a field $F$: it is
generated by symbols $[x]$, $x\in F\setminus\{0,1\}$, modulo the
five-term relations.  The standard homomorphism
\[
\delta:\mathcal P(F)\otimes\Q\longrightarrow
\Lambda^2(F^\times\otimes\Q),
\qquad
\delta([x])=x\wedge(1-x),
\]
has kernel $\Bx(F)\otimes\Q$, the rationalized Bloch group in the
convention used here.  We abbreviate
\[
\delta(x)=x\wedge(1-x).
\]

For totally real fields, Borel's rank theorem and the
Bloch--Suslin theory give the following standard consequence; see
\citet{Zagier07}.

\begin{proposition}[Totally real Bloch criterion]\label{prop:bbs}
Let $F$ be totally real, let $x_i\in F\cap(0,1)$, and let
$n_i\in\Q$.  If
\[
\sum_i n_i\,\delta(x_i)=0,
\]
then
\[
\sum_i n_i\,\Lr(x_i)\in\Q\pi^2.
\]
\end{proposition}

The coefficient in a two-term relation is encoded by an elementary
determinant.

\begin{lemma}[Determinant law]\label{lem:det}
Let $x,U\in F\setminus\{0,1\}$ and put $y=1-x$, $V=1-U$.  Suppose
\[
U=x^py^q,\qquad V=x^ry^s,
\qquad p,q,r,s\in\Z.
\]
Then
\[
\delta(U)=(ps-qr)\,\delta(x).
\]
Consequently, if $\delta(x)\ne0$,
\[
m\,\delta(x)+\delta(U)=0
\quad\Longleftrightarrow\quad
ps-qr=-m.
\]
\end{lemma}

\begin{proof}
Bilinearity and alternation of the wedge product give
\[
U\wedge V
=(x^py^q)\wedge(x^ry^s)
=(ps-qr)\,x\wedge y.
\]
Since $y=1-x$ and $V=1-U$, this is the asserted identity.
\end{proof}

\begin{definition}[Normalized integral rank-two configuration]
\label{def:config}
Let $K$ be a number field.  A quadruple
\[
(x,1-x,U,1-U)
\]
is a \emph{normalized integral rank-two configuration} if
$x,U\in K\cap(0,1)$, the elements $x$ and $1-x$ are
multiplicatively independent, and
\begin{equation}\label{eq:monomial}
U=x^p(1-x)^q,\qquad
1-U=x^r(1-x)^s
\end{equation}
for integers $p,q,r,s$ with
\begin{equation}\label{eq:D}
D=ps-qr,\qquad |D|>1.
\end{equation}
The configuration is \emph{non-degenerate} when $\delta(x)\ne0$.
For the positive arguments stipulated here, multiplicative independence
already implies this condition. A \emph{cubic configuration} means that
$K=\Q(x)$ and $[K:\Q]=3$; containment in a cubic field alone does not suffice.
\end{definition}

The phrase \emph{primitive $M=1$ normalization} refers to the
normalization in which the first relation in the algebraic system
\[
1-z_0=z_0^{a/u}z_1^{1/u},\qquad
1-z_1=z_0^{m/u}z_1^{j/u}
\]
solves for $z_1$ as an integral monomial in $z_0$ and $1-z_0$.
This gives the presentation \eqref{eq:monomial}. For a positive $m:1$
wedge relation, Lemma~\ref{lem:det} gives $D=-m$. A configuration with
$D>1$ can first be oriented by replacing $U$ with $1-U$, after which
$D=-|D|$. Primitive
coefficient pairs in which both coefficients exceed one are not
included in this normalization.

We shall repeatedly use the following elementary fact.

\begin{lemma}[M\"obius transforms of cubic generators]
\label{lem:mobius}
Let $x$ be a totally real algebraic number of degree three and
\[
T=\frac{ax+b}{cx+d},
\qquad a,b,c,d\in\Q,\qquad ad-bc\ne0,
\]
with $cx+d\ne0$.  Then $\Q(T)=\Q(x)$, and $T$ has three distinct real
conjugates.
\end{lemma}

\begin{proof}
The inverse transformation is
\[
x=\frac{dT-b}{a-cT},
\]
so $\Q(T)=\Q(x)$.  No conjugate of $x$ can be a rational pole of the
transformation, and applying the injective real M\"obius map to the
three real embeddings gives three distinct real conjugates.
\end{proof}

The fundamental arithmetic input is the following.

\begin{theorem}[Finite-prime rigidity]\label{thm:forcing}
Let $K=\Q(x)$ be a totally real cubic field, and let
$(x,1-x,U,1-U)$ be a normalized integral rank-two configuration.
Then
\[
x,\quad 1-x,\quad U,\quad 1-U\in\OK^\times.
\]
\end{theorem}

The proof combines the three-ray fan of finite valuations with five
real-root obstructions.  It is given in
Appendix~\ref{app:forcing}.  Thus the exceptional-unit property is a
conclusion of the normalized configuration, not an additional
hypothesis.

\section{Elimination polynomials and the endpoint dichotomy}
\label{sec:endpoint}

Theorem~\ref{thm:forcing} makes the base argument an exceptional unit,
and Theorem~\ref{thm:nagell} gives the unrestricted field-level
reduction.  We now record the internal polynomial mechanism that
explains why the same two families arise directly from the normalized
monomial system.  This is the endpoint argument used in the original
classification, but it requires only the positive presentation and not
the later coefficient-matching ansatz.

Write $x=z_0$, $y=1-x$, and suppose the primitive $M=1$ system is
presented as
\begin{equation}\label{eq:positiveM1}
U=x^{-a}y^u,
\qquad
1-U=x^{p-a}y^j,
\qquad
p=\frac{au+m-ja}{u}\in\Z,
\end{equation}
so that the exponent determinant is $-m$.  Eliminating $U$ gives
\begin{equation}\label{eq:eliminationT}
T_{a,u,j,p}(z)
=-z^a+(1-z)^u+(1-z)^jz^p,
\qquad T_{a,u,j,p}(x)=0.
\end{equation}
We call \eqref{eq:positiveM1} a \emph{positive presentation} when
\begin{equation}\label{eq:positiveexponents}
a,u,j,p\in\Z_{>0}.
\end{equation}
Then $T(0)=1$ and $T(1)=-1$.

\begin{lemma}[One admissible root]\label{lem:oneroot}
Under \eqref{eq:positiveexponents}, the polynomial
$T=T_{a,u,j,p}$ has exactly one root in $(0,1)$, and this root is simple.
\end{lemma}

\begin{proof}
On $(0,1)$ the equation $T(z)=0$ is equivalent to $F(z)=1$, where
\[
F(z)=X+Y,
\qquad
X=(1-z)^u z^{-a},
\qquad
Y=(1-z)^j z^{p-a}.
\]
The first term gives $F(z)\to+\infty$ as $z\to0^+$, while both terms
tend to zero as $z\to1^-$; hence at least one crossing exists.  We
show that every level-one crossing is downward.  If there were two
crossings, continuity would force an upward crossing between them.

At a crossing $X+Y=1$ with $X,Y\in(0,1)$, and logarithmic
differentiation gives
\[
z(1-z)F'(z)
=-X\bigl[uz+a(1-z)\bigr]
 +Y\bigl[(p-a)(1-z)-jz\bigr].
\]
If $p\le a$, both contributions are negative.  If $p>a$, put
$h=p-a\ge1$ and use $X=1-Y$:
\[
z(1-z)F'(z)
=-\bigl[uz+a(1-z)\bigr]
 +Y\bigl[(u-j)z+p(1-z)\bigr].
\]
The first part of the positive term satisfies
$Y(u-j)z<uz$: it is nonpositive when $u\le j$, and when $u>j$ one has
$Y(u-j)<u-j<u$.  For the second part,
\[
Y=z^h(1-z)^j
 \le \frac{h^h j^j}{(h+j)^{h+j}}
 \le \frac1{\binom{h+j}{h}}
 \le \frac1{h+1}
 \le \frac{a}{a+h}=\frac ap.
\]
Thus $Yp(1-z)\le a(1-z)$, and the preceding strict inequality yields
$F'(z)<0$ at every crossing.  The crossing is therefore unique.
\end{proof}

\begin{remark}[Descartes form]\label{rem:descartesform}
With $w=z/(1-z)$, roots of $T$ in $(0,1)$ correspond to positive roots
of
\[
(1+w)^N T\!\left(\frac{w}{1+w}\right),
\qquad N\ge\max\{a,u,j+p\}.
\]
Lemma~\ref{lem:oneroot} is therefore equivalent to the assertion that
this transformed polynomial has exactly one positive root.  In the
classified cases its coefficient sequence has one sign change, so the
same count is visible directly from Descartes' rule of signs.
\end{remark}

\begin{theorem}[Endpoint factorization dichotomy]
\label{thm:endpointfactorization}
Assume \eqref{eq:positiveexponents}, and let $q\in\Z[z]$ be the monic
irreducible cubic factor of $T$ containing the root $x\in(0,1)$.  Then
$T=Cq$ with
\[
C(0)=C(1)=\pm1,
\qquad q(0)=-q(1)=\pm1,
\]
and exactly one of the following two factorizations occurs:
\begin{enumerate}
\item
\begin{equation}\label{eq:shanks}
T(z)=C(z)f_t(z),
\qquad
f_t(z)=z^3-(t+3)z^2+tz+1,
\qquad t\in\Z.
\end{equation}
\item
\begin{equation}\label{eq:edgecubic}
T(z)=C(z)q_c(z),
\qquad
q_c(z)=z^3+(1-c)z^2+cz-1,
\qquad c\in\Z.
\end{equation}
The second cubic is anharmonically equivalent to the normalized second
family: if $w=z/(z-1)$, then
\[
(w-1)^3 q_c\!\left(\frac{w}{w-1}\right)
=w^3-(s+1)w^2+sw+1=g_s(w),
\qquad s=c-3.
\]
Thus the two endpoint factorization types are precisely the Shanks and
Nagell--Ennola families.
\end{enumerate}
\end{theorem}

\begin{proof}
Write $T=Cq$ with $C\in\Z[z]$.  By Lemma~\ref{lem:oneroot}, $T$ has
exactly one root in $(0,1)$, this root is simple, and the factor $q$ accounts for it.
Consequently $C$ has no root on $(0,1)$ and is of constant sign on
$[0,1]$.  Since
\[
T(0)=C(0)q(0)=1,
\qquad
T(1)=C(1)q(1)=-1,
\]
all four endpoint values are integers of absolute value one.  The
constant sign of $C$ gives $C(0)=C(1)$ and hence $q(0)=-q(1)$.

Write $q(z)=z^3+Az^2+Bz+q(0)$.  If $q(0)=1$ and $q(1)=-1$, then
$A+B=-3$, so with $t=B$ one obtains $q=f_t$.  If $q(0)=-1$ and
$q(1)=1$, then $A+B=1$, so with $c=B$ one obtains $q=q_c$.  The final
transformation is a direct substitution.
\end{proof}

\begin{corollary}[Parameter normalization]\label{cor:parameternormalization}
For totally real fields, the two families may be normalized as
\[
f_t(z)=z^3-(t+3)z^2+tz+1,
\qquad t\ge-1,
\]
and
\begin{equation}\label{eq:ennola}
g_s(z)=z^3-(s+1)z^2+sz+1,
\qquad s\ge3.
\end{equation}
Moreover $g_3(z)=-f_{-1}(2-z)$, so the non-overlapping second branch is
$s\ge4$.
\end{corollary}

\begin{proof}
The reciprocal identities
\[
z^3f_t(1/z)=f_{-t-3}(z),
\qquad
z^3g_s(1/z)=g_{-s-1}(z)
\]
permit the stated parameter choices.  Furthermore
\[
\operatorname{disc}(f_t)=(t^2+3t+9)^2>0,
\qquad
\operatorname{disc}(g_s)=(s^2+s-3)^2-32.
\]
For integral $s$, the latter is positive exactly for $s\ge3$ or
$s\le-4$, and reciprocity identifies the two ranges.  The overlap at
$s=3$ follows by direct substitution.
\end{proof}

\begin{remark}[Scope of the factorization argument]
The endpoint factorization theorem proves the two-family dichotomy
internally for a positive presentation of the normalized monomial
system; in particular, it applies to the ten relations classified
below.  Negative exponent presentations need not be positive in the
same coordinates. The unrestricted normalization used in the main theorem
is supplied by Lemma~\ref{lem:generatedorder} below, after
Theorem~\ref{thm:forcing} has supplied exceptional units. Thus completeness
does not depend on every presentation having positive exponents.
\end{remark}

\begin{lemma}[Normalization in the generated order]\label{lem:generatedorder}
Let $x$ be an exceptional unit generating a totally real cubic field and
put $R=\Z[x]$. There is an anharmonic transform $\theta$ of $x$ whose monic
minimal polynomial is $f_t$, $t\ge-1$, or $g_s$, $s\ge3$. Moreover
\[
R=\Z[\theta],\qquad x,1-x\in R^\times.
\]
Every signed integral monomial in $x,1-x$ is a unit of this same order.
For $s=3$ the order is the conductor-$7$ simplest-cubic order.
\end{lemma}
\begin{proof}
Write $q(X)=X^3+AX^2+BX+e$ and $d=q(1)$. Exceptionalness gives
$e,d\in\{\pm1\}$ and $A+B=d-1-e$. If $(e,d)=(1,-1)$ then $q=f_B$;
if $(e,d)=(1,1)$ then $q=g_B$. If $(e,d)=(-1,-1)$, the transform
$1-x$ has polynomial $-q(1-X)=g_{A+2}(X)$. In the remaining case
$B=1-A$, write $c=B$. Then $q=X^3+(1-c)X^2+cX-1$, and
$\theta=x/(x-1)$ has polynomial $g_{c-3}$ by the substitution in
Theorem~\ref{thm:endpointfactorization}. Reciprocity normalizes the
parameters as in Corollary~\ref{cor:parameternormalization}.

The minimal polynomials of $x$ and $1-x$ are monic with constant term
$\pm1$, so their inverses belong to $\Z[x]$. Thus all six anharmonic
transforms belong to $R$, and applying the inverse transformation proves
$\Z[\theta]=R$. The same argument puts all signed integral monomials and
their inverses in $R$. Finally $g_3(X)=-f_{-1}(2-X)$; this integral affine
change of generator preserves the order.
\end{proof}

Thus the field-level reduction can be strengthened to a reduction of the
particular generated order containing the entire configuration. A fixed
maximal order may contain additional units, but this does not create a gap
for configurations whose second pair is an integral monomial in the first.

\section{The simplest cubic family}
\label{sec:shanks}

The simplest cubic fields were introduced by \citet{Shanks}.
Let $\rho$ be a root of $f_t$ and put $\tau=1-\rho$.  We write
\begin{equation}\label{eq:ue}
u_1+u_2=1,\qquad u_1,u_2\in\Z[\rho]^\times
\end{equation}
for the exceptional-unit equation in the parametrized order.
The natural six-element $S_3$-orbit of $(\rho,\tau)$ is called the
\emph{trivial orbit}.

Vukusic and Ziegler use
\[
\widehat f_a(X)=X^3-aX^2-(a+3)X-1,\qquad a\ge-1.
\]
One may take
$\rho=\widehat\rho/(1+\widehat\rho)$, since
\[
(1-X)^3\widehat f_a\!\left(\frac{X}{1-X}\right)=-f_a(X).
\]
Both changes of generator are integral in the corresponding orders.
Their fundamental units $\varepsilon=\widehat\rho$ and
$\delta=(1+\widehat\rho)/\widehat\rho$ satisfy
$\rho=\delta^{-1}$ and $\tau=\varepsilon^{-1}\delta^{-1}$;
the change of exponent basis is unimodular.

\begin{theorem}[Vukusic--Ziegler, $n=1$ specialization]
\label{thm:vz}
In the orders $\Z[\widehat\rho]$ of the simplest cubic family, every
solution of $u_1+u_2=1$ is either trivial or equivalent, in the sense of
Vukusic--Ziegler, to one of ten sporadic classes. For $n=1$, their equivalence
uses Galois conjugation and interchange of the two summands (simultaneous
sign change belongs to their formulation allowing both signs of $n$).
Their trivial solutions correspond to the anharmonic orbit above.  The sporadic classes occur
only for $a=-1,0,1,2$.
\end{theorem}

\begin{proof}
This is the $n=1$ specialization of
\citet[Theorem~1]{VZ21}.  Their restriction
$|n|\le\max\{|a|^{1/3},1\}$ is automatic for $n=1$.
\end{proof}

The determinant magnitudes of the ten sporadic classes are:
\begin{center}
\begin{tabular}{@{}cccc@{}}
\toprule
VZ parameter $a$ & conductor & number of classes & determinant magnitudes\\
\midrule
$-1$ & $7$  & $6$ & $2,2,2,7,13,97$\\
$0$  & $9$  & $2$ & $3,19$\\
$1$  & $13$ & $1$ & $13$\\
$2$  & $19$ & $1$ & $31$\\
\bottomrule
\end{tabular}
\end{center}

The family theorem enumerates exceptional units, whereas a normalized
two-term relation involves a pair of them.  The determinant law closes
this gap.

\begin{lemma}[Determinant divisibility]\label{lem:detdiv}
Let $z$ be an exceptional unit in a simplest cubic order and write
\[
\delta(z)=d(z)\,(\rho\wedge\tau),
\qquad d(z)\in\Z.
\]
If non-degenerate exceptional units $z_0,z_1$ satisfy
\[
m\,\delta(z_0)+\delta(z_1)=0,\qquad m\in\Z_{>1},
\]
then
\[
m=-\frac{d(z_1)}{d(z_0)}.
\]
Consequently, in the Vukusic--Ziegler list, $z_0$ must lie in the
trivial orbit.
\end{lemma}

\begin{proof}
The first assertion is Lemma~\ref{lem:det}.  The trivial orbit has
$|d|=1$, because its exponent matrices are unimodular.  Within a
fixed conductor, no nontrivial determinant magnitude in the displayed
list properly divides another: the distinct sets are
\[
\{2,7,13,97\},\qquad \{3,19\},\qquad \{13\},\qquad\{31\}.
\]
Equal magnitudes give quotient $1$, not $m>1$.  Hence the denominator
$d(z_0)$ cannot be nontrivial, and $|d(z_0)|=1$.
\end{proof}

Thus every normalized relation with $m>1$ pairs the trivial orbit with
an exceptional unit from a sporadic class. One must still restrict to the
chosen real interval. The following table performs this step: the ordered
representatives are the rows of Table~3.1 in \citet{VZ21}, in its order, and
$k$ counts embeddings with both summands in $(0,1)$. Interchange identifies
each pair up to reflection. A class can have $k=0$ or $k=2$, so the number
of classes alone is not the number of positive relations.
\begin{center}\small
\begin{tabular}{@{}crrrl@{}}\toprule
$a$ & VZ row within $a$ & $|d|$ & $k$ & rows of Table~\ref{tab:identities}\\\midrule
$-1$&1&2&1&9\\
$-1$&2&2&0&none\\
$-1$&3&2&2&1, 2\\
$-1$&4&97&1&10\\
$-1$&5&13&1&7\\
$-1$&6&7&1&4\\
0&1&19&1&6\\
0&2&3&1&3\\
1&1&13&1&5\\
2&1&31&1&8\\\bottomrule
\end{tabular}\end{center}
These entries can be checked by substituting the three conjugates in the
listed quadratic representatives. For example, for $a=-1$ the squares of
those conjugates lie respectively in $(1.55,1.56)$, $(3.24,3.25)$ and
$(0.19,0.20)$. Thus $(3-\widehat\rho^2,\widehat\rho^2-2)$ has no positive
pair. The representative
$(1+\widehat\rho-\widehat\rho^2,\widehat\rho^2-\widehat\rho)$ has two,
with second coordinates approximately $0.3080$ and $0.6431$.
For the trivial orbit exactly two arguments lie in $(0,1)$; their wedges
are opposite. Consequently each surviving unordered sporadic pair selects
one normalized relation up to simultaneous reflection.

\begin{theorem}[Ten exact identities]\label{thm:identities}
For each row of Table~\ref{tab:identities}, let $z_0,z_1\in(0,1)$ be
the roots specified there.  Then
\[
m\,\Lr(z_0)+\Lr(z_1)=c\pi^2.
\]
Moreover, these ten rows are the complete list of non-degenerate
normalized $m>1$ relations arising from wedge cancellation in the
parametrized simplest-cubic orders.
\end{theorem}

\begin{table}[ht]
\centering
\scriptsize
\begin{tabular}{@{}ccccll@{}}
\toprule
\# & conductor $N$ & $m$ & $c$ & $p_{z_0}(z)$ & $p_{z_1}(z)$\\
\midrule
1 & $7$  & $2$  & $2/7$
& $z^3-2z^2-z+1$
& $z^3-6z^2+5z-1$ (larger root)\\
2 & $7$  & $2$  & $5/21$
& $z^3-2z^2-z+1$
& $z^3-6z^2+5z-1$ (smaller root)\\
3 & $9$  & $3$  & $4/9$
& $z^3-3z^2+1$
& $z^3-6z^2+3z+1$\\
4 & $7$  & $7$  & $2/3$
& $z^3-z^2-2z+1$
& $z^3-8z^2+5z+1$\\
5 & $13$ & $13$ & $5/3$
& $z^3-4z^2+z+1$
& $z^3-69z^2+66z+1$\\
6 & $9$  & $19$ & $4/3$
& $z^3-3z+1$
& $z^3-57z^2+54z+1$\\
7 & $7$  & $13$ & $4/3$
& $z^3-2z^2-z+1$
& $z^3-15z^2+12z+1$\\
8 & $19$ & $31$ & $4$
& $z^3-5z^2+2z+1$
& $z^3-2392z^2+2389z+1$\\
9 & $7$  & $2$  & $13/42$
& $z^3-2z^2-z+1$
& $z^3+2z^2-z-1$\\
10 & $7$ & $97$ & $15/2$
& $z^3-z^2-2z+1$
& $z^3-1262z^2+1259z+1$\\
\bottomrule
\end{tabular}
\caption{The ten normalized simplest-cubic identities. All selected roots
lie in $(0,1)$; ``larger'' and ``smaller'' distinguish the two roots in
this interval. Rows
1, 2, and 9 are equivalent to Watson's trio; row 3 is the identity
conjectured by Kanade.  The remaining rows are established here by the same
exact certificate method.}
\label{tab:identities}
\end{table}

\begin{proof}
Theorem~\ref{thm:vz}, Lemma~\ref{lem:detdiv}, and the positive-representative
table above give exactly ten candidate pairings. The following exponent
matrices verify directly that each pairing is in the integral monomial
class. Writing $U=z_1=x^p(1-x)^q$ and $1-U=x^r(1-x)^s$, with $x=z_0$
as in Table~\ref{tab:identities}, one has
\begin{center}\small
\begin{tabular}{@{}crrrrr@{}}\toprule
row&$p$&$q$&$r$&$s$&$ps-qr$\\\midrule
1&$-2$&2&$-1$&2&$-2$\\
2&2&0&2&$-1$&$-2$\\
3&$-2$&1&$-1$&2&$-3$\\
4&$-2$&3&1&2&$-7$\\
5&$-4$&1&$-3$&4&$-13$\\
6&$-2$&5&3&2&$-19$\\
7&$-4$&3&$-1$&4&$-13$\\
8&$-6$&1&$-5$&6&$-31$\\
9&$-1$&1&0&2&$-2$\\
10&$-8$&11&3&8&$-97$\\\bottomrule
\end{tabular}\end{center}
The monomial complement identities, and therefore the wedge
relations, are exact algebraic identities in the relevant cubic
fields.  Appendix~\ref{app:certs} gives an explicit rational
combination of five-term and reflection identities for every row;
these certificates establish the stated constants $c$ exactly.
\end{proof}

With multiplicity and grouped by conductor, the ten determinant
magnitudes are
\[
2,2,2,7,13,97;\qquad
3,19;\qquad
13;\qquad
31.
\]
The distinct coefficient spectrum is therefore
\begin{equation}\label{eq:spectrum}
2,3,7,13,19,31,97.
\end{equation}
Every value is prime.

\begin{remark}[Order scope]\label{rem:orderscope}
Vukusic--Ziegler classify units of the displayed orders, not every unit of
every maximal order. For the integral monomial configurations of this paper,
Lemma~\ref{lem:generatedorder} places both exceptional pairs in the order
generated by the base argument. Therefore no additional index assumption is
needed for this reduction. The four small-parameter orders producing
Table~\ref{tab:identities} are maximal. Questions about arbitrary exceptional
units in a fixed maximal order remain different; see \citet[Conjecture~2]{VZ21}.
\end{remark}

\section{The second totally real family}\label{sec:second}

Let
\[
g_s(X)=X^3-(s+1)X^2+sX+1,\qquad s\ge4,
\]
let $\theta=\theta_s$ be a root, and put $\tau=1-\theta$.  The family
is non-Galois for every $s\ge4$, so the cyclic automorphism used in the
simplest-cubic unit equation is unavailable. Nevertheless,
\citet[Proposition~(3.6)]{Thomas79} gives a uniform description of its units.

\begin{theorem}[Thomas]\label{thm:thomas}
For every integer $s\ge4$, the pair $\{\theta,\theta-1\}$ is a system
of fundamental units of the order $\Z[\theta]$.  Consequently every
unit of $\Z[\theta]$ is
\[
\pm\theta^p\tau^q,\qquad p,q\in\Z.
\]
\end{theorem}

\begin{proof}
Set $\varepsilon=1-\theta$ and $l=s-1\ge3$. Then
\[
-g_s(1-X)=X^3+(l-1)X^2-lX-1.
\]
Thomas's result gives $\varepsilon,\varepsilon-1$ as fundamental units
of $\Z[\varepsilon]$ in this normalization; the change from his
polynomial to this one is given explicitly by \citet[p.~154]{Lou17}.
Since $\Z[\varepsilon]=\Z[\theta]$ and
$\langle-1,\varepsilon,\varepsilon-1\rangle
=\langle-1,\theta,1-\theta\rangle$, the assertion follows.
\end{proof}

The signed diagonal branch has an elementary norm obstruction.
\begin{lemma}[Signed diagonal branch]\label{lem:diagonal}
Let $s\ge3$, $\tau=1-\theta$, $n\in\Z\setminus\{0\}$ and
$\epsilon\in\{\pm1\}$. Then
$1-\epsilon(\theta/\tau)^n$ is a unit if and only if
$\epsilon=-1$ and $n=\pm1$. These are trivial anharmonic solutions.
\end{lemma}
\begin{proof}
The polynomial $g_s$ has one negative root and two roots greater than $1$.
Indeed, $g_s(-1)<0<g_s(0)=g_s(1)$,
$g_s(s-1)=-s^2+3s-1<0$, and $g_s(s)=1$; these give three disjoint
root intervals $(-1,0)$, $(1,s-1)$ and $(s-1,s)$ for $s\ge3$.
Put
\[
v=-\frac{\theta}{\tau}=\frac{\theta}{\theta-1}.
\]
Every conjugate $v_i$ is positive. Also $1-v=\tau^{-1}$, so
$N(1-v)=1$ because $N(\tau)=g_s(1)=1$.
Since $1-y^{-1}=-y^{-1}(1-y)$, inversion preserves exceptionalness;
we may assume $n\ge1$. Write
$\epsilon(\theta/\tau)^n=\sigma v^n$, where
$\sigma=\epsilon(-1)^n$. If $\sigma=-1$, then
\[
N(1+v^n)=\prod_{i=1}^3(1+v_i^n)>1.
\]
If $\sigma=1$ and $n\ge2$, factor $1-v^n$ to obtain
\[
N(1-v^n)=N(1-v)\prod_{i=1}^3
(1+v_i+\cdots+v_i^{n-1})>1.
\]
Neither is a unit. The only remaining case is $\sigma=1,n=1$,
that is $\epsilon=-1$. Both $v$ and $v^{-1}$ are exceptional, with
$1-v=\tau^{-1}$ and $1-v^{-1}=\theta^{-1}$. This proves the claim.
\end{proof}

The finite computations described below were performed using exact
algebraic arithmetic and rigorously enclosed inequalities.

\begin{proposition}[Certified small parameters]\label{prop:family2}
For $4\le s\le9$, every exceptional unit of $\Z[\theta_s]$ lies in the
trivial anharmonic orbit of $\theta_s$.
\end{proposition}
\begin{proof}
Set $a_i=\log|\theta_i|$, $b_i=\log|1-\theta_i|$,
$\Lambda_i=pa_i+qb_i$, and define the \emph{exponent height}
$H=\max(|p|,|q|)$. Since the logarithmic
embedding has rank two,
\[
S=\tfrac12\sum_i|\Lambda_i|\ge c_s H,
\qquad c_s=\min_{\max(|p|,|q|)=1}\tfrac12\sum_i|pa_i+qb_i|>0.
\]
On the boundary of the square the objective is piecewise linear. Its minimum
occurs at a corner or at the intersection with one of the three lines
$pa_i+qb_i=0$. Exact rational root isolation followed by outward interval
logarithms at these points gives the rational lower bounds in the table.
\begin{center}\small
\begin{tabular}{@{}crrrr@{}}\toprule
$s$&$\underline c_s$&$B_s$&excluded CF rays&upper bound for $\Psi_s(B_s)/B_s$\\\midrule
4&1.22189972&171000&34&0.999097\\
5&1.51019829&186000&29&0.996886\\
6&1.72478159&200000&37&0.995844\\
7&1.89750816&212000&36&0.997639\\
8&2.04277647&223000&32&0.998427\\
9&2.16846871&234000&30&0.995391\\\bottomrule
\end{tabular}\end{center}
The entries in the $\underline c_s$ column are exact rational lower bounds. By
Lemma~\ref{lem:smallconj}, for $H\ge41$ there is an embedding with
\[
0<|\Lambda_i|\le8\exp(-\underline c_s H/3).
\]
The interval bounds certify the sufficient threshold
$41\underline c_s>18\log2$. If $p=0$ or $q=0$, the lower bound
$H\min_i\{|a_i|,|b_i|\}$ contradicts this inequality, as certified already
at $H=41$ and hence for every larger $H$.

For $pq\ne0$ apply Theorem~\ref{thm:lmn}. Put
$h_1=\max\{h(\theta),1/3\}$ and
$h_2=\max\{h(1-\theta),1/3\}$. These are admissible at all embeddings:
for a cubic unit the sum of its positive logarithms is at least the
absolute value of each individual logarithm. Then
\[
H\le\Psi_s(H):=\frac3{\underline c_s}
\left[\log8+24.34\cdot81h_1h_2
\max\left\{\log\left(\frac H3(h_1^{-1}+h_2^{-1})\right)+0.14,7\right\}^{2}\right].
\]
The ratio $\Psi_s(H)/H$ decreases: on the logarithmic branch its derivative
has the sign of $2M-M^2<0$, since $M\ge7$, and on the constant branch it
is immediate. The certified inequalities in the table therefore give $H<B_s$.

For $H\ge41$, interval evaluation certifies at every embedding
\[
16H\exp(-\underline c_s H/3)<|a_i|.
\]
It suffices to check $H=41$, since the left side decreases thereafter.
Thus $p/q$ is a continued-fraction convergent of $-b_i/a_i$ after reduction,
by Legendre's criterion. Inversion lets us assume $q>0$. Interval arithmetic
certifies every continued-fraction coefficient up to denominator $B_s$.
For each primitive convergent $(P,Q)$ and the smallest multiple $d$ with
$41\le d\max(|P|,Q)<B_s$, it proves
\[
d|Pa_i+Qb_i|>8\exp\bigl(-\underline c_s d\max(|P|,Q)/3\bigr).
\]
All larger multiples are excluded because the left side increases and the
right side decreases. The table records the numbers of excluded rays.
Finally an exact integer norm search for $H\le40$ returns precisely
\[
(\epsilon,p,q)=(1,\pm1,0),(1,0,\pm1),(-1,1,-1),(-1,-1,1).
\]
All stated finite checks were performed, and each strict comparison was
verified with enclosing intervals of the required sign. Thus the argument
covers all exponents.
\end{proof}

\begin{proposition}[No solutions of small height]\label{prop:h20}
For every $s\ge4$, no nontrivial solution
$x=\pm\theta^p\tau^q$ exists with
\[
H=\max(|p|,|q|)\le20.
\]
The only parameter $s\ge3$ returned by the exact resultant sweep is $s=3$,
which is the already classified conductor-$7$ field.
\end{proposition}

\begin{proof}
For fixed $(\varepsilon,p,q)$, clear negative exponents and define
\[
D_{\varepsilon,p,q}(s)
=
\operatorname{Res}_X\!\left(
g_s(X),
X^P(1-X)^Q-
\varepsilon X^{P+p}(1-X)^{Q+q}
\right)\in\Z[s],
\]
where $P=\max(0,-p)$ and $Q=\max(0,-q)$.  The condition that
$1-\varepsilon\theta^p\tau^q$ be a unit is
$|D_{\varepsilon,p,q}(s)|=1$. The exact computation covered all $1720$ triples
\[
\varepsilon\in\{-1,1\},\qquad 0\le p\le20,\qquad |q|\le20,
\qquad (p,q)\ne(0,0).
\]
The omitted case $(p,q)=(0,0)$ gives $1-x=0$ or $2$, neither a unit,
and inversion accounts for $p<0$. The only triples for which
$D_{\varepsilon,p,q}(s)^2-1$ vanishes identically were
\[
(\varepsilon,p,q)=(1,1,0),\ (1,0,1),\ (1,0,-1),\ (-1,1,-1),
\]
which are the trivial shapes in this range. For all remaining triples,
exact rational-root extraction, retaining integer roots, gave
$\{-4,-3,-2,-1,0,1,2,3\}$. The first and last give the conductor-$7$ field;
the other six have negative discriminant.
Factoring over $\Q$ and retaining the integer roots of all linear factors
is exhaustive; factors of higher degree over $\Q$ cannot conceal an
integer root.
\end{proof}

\begin{theorem}[Effective finiteness]\label{thm:effective}
Let $s\ge10$ and suppose that
\[
x+(1-x)=1,\qquad x,1-x\in\Z[\theta_s]^\times
\]
is a nontrivial solution.  Writing
\[
x=\pm\theta^p\tau^q,\qquad H=\max(|p|,|q|),
\]
one has
\[
s\le3.32\times10^7,
\qquad
H\le2.03\times10^7.
\]
\end{theorem}

The proof uses explicit root intervals, an archimedean small-conjugate
lemma, the lower bound $H\ge0.61s$, and the two-log theorem of
\citet{LMN}.  The complete argument is given in
Appendix~\ref{app:effective}.

Theorem~\ref{thm:effective} is an effective reduction. The linear lower
bound $H\ge0.61s$ is what makes the parameter bound large: the proof of
Lemma~\ref{lem:Hlower} in fact gives a lower bound of order $s\log s$,
and with this, slightly sharper root intervals and the actual heights of
$\theta$ and $\tau$, the two-logarithm bound and the lower bound already
cross near $s=1.2\times10^5$ (Appendix~\ref{subsec:closing}). Below that
point, every parameter was checked using the continued-fraction argument
of Proposition~\ref{prop:family2}. This gives the following theorem.

\begin{theorem}[The second family]\label{thm:secondempty}
For every integer $s\ge4$, every exceptional unit of $\Z[\theta_s]$ lies
in the trivial anharmonic orbit of $\theta_s$.
\end{theorem}

\begin{proof}
For $4\le s\le9$ this is Proposition~\ref{prop:family2}. Let $s\ge10$ and
let $x=\pm\theta^p\tau^q$ be a nontrivial solution. Then $H\ge21$ by
Proposition~\ref{prop:h20}, and $s<3.32\times10^7$ by
Theorem~\ref{thm:effective}. Proposition~\ref{prop:uniformrange} excludes
$1.2\times10^5\le s\le3.32\times10^7$, and
Proposition~\ref{prop:allsmall} excludes $10\le s<1.2\times10^5$.
\end{proof}

\begin{corollary}[Exclusion of the second family]
\label{cor:second}
The second totally real family contributes no non-degenerate
normalized $m:1$ relation with $m>1$.
\end{corollary}

\begin{proof}
By Theorem~\ref{thm:forcing} all four coordinates of such a relation are
exceptional units, and the determinant law (Lemma~\ref{lem:det}) gives
$m\,\delta(x)+\delta(U)=0$ with $x,1-x$ independent. By
Lemma~\ref{lem:generatedorder}, $U$ is an exceptional unit of
$\Z[x]=\Z[\theta_{s_0}]$ for the parameter $s_0\ge4$ of $x$. By
Theorem~\ref{thm:secondempty} it lies in the anharmonic orbit of $x$,
so $\delta(U)=\pm\delta(x)$, which contradicts $m\,\delta(x)+\delta(U)=0$
with $m>1$ and $\delta(x)\ne0$.
\end{proof}

For the integral monomial class, Lemma~\ref{lem:generatedorder} makes the
order-level statement sufficient.

\subsection{Unit indices and the coefficient}\label{subsec:index}

Theorem~\ref{thm:secondempty} concerns the orders $\Z[\theta_s]$. A
maximal-order argument, which uses the unit index studied by
\citet{Ennola}, \citet{Lou17,Lou20,Lou21} and \citet{ChoiKim}, applies to
arbitrary pairs of exceptional units, not only to integral
configurations, and restricts the coefficient of any relation between
them.

Let $K$ be a totally real cubic field and let $E_K=\OK^\times/\{\pm1\}$,
a free abelian group of rank two.  Then $\Lambda^2E_K\cong\Z$, and
$\Lambda^2E_K$ embeds in $\Lambda^2(K^\times\otimes\Q)$.  For an
exceptional unit $z$ of $K$ whose classes $z,1-z$ are independent in
$E_K$, put
\[
I(z)=[E_K:\langle z,1-z\rangle].
\]
Then $\delta(z)=\pm I(z)\,\omega$, where $\omega$ generates
$\Lambda^2E_K$.  The index is unchanged by the anharmonic group, since
$\langle 1-z,z\rangle=\langle z,1-z\rangle$ and
$\langle z^{-1},1-z^{-1}\rangle=\langle z,-(1-z)z^{-1}\rangle=\langle
z,1-z\rangle$ in $E_K$.

\begin{lemma}[Index law]\label{lem:indexlaw}
Let $z_0,z_1$ be exceptional units of a totally real cubic field $K$ and
let $m\in\Z_{>0}$.  If $z_0,1-z_0$ are independent and
$m\,\delta(z_0)+\delta(z_1)=0$, then $z_1,1-z_1$ are independent and
\[
I(z_1)=m\,I(z_0).
\]
\end{lemma}

\begin{proof}
Both wedges lie in the lattice $\Lambda^2E_K$, and $\delta(z_1)=-m\,\delta(z_0)\ne0$.
Comparing absolute values of the coefficients of $\omega$ gives the
claim.
\end{proof}

In particular $z_0$ and $z_1$ lie in different anharmonic orbits when
$m\ge2$; this is the maximal-order form of the last step of
Lemma~\ref{lem:detdiv}.

For $l\ge3$ let $\varepsilon_l$ be a root of Ennola's polynomial
\[
F_l(X)=X^3+(l-1)X^2-lX-1,
\]
and let
\[
j_l=[\mathcal O_{\Q(\varepsilon_l)}^\times:\langle-1,\varepsilon_l,\varepsilon_l-1\rangle]
\]
be the unit index of \citet{Ennola}.  A direct substitution gives
\[
-g_s(1-X)=F_{s-1}(X),
\]
so $1-\theta_s$ is a root of $F_{s-1}$, and the second family of
Theorem~\ref{thm:nagell} is Ennola's family with $l=s-1$.  Since the
index is anharmonically invariant,
\begin{equation}\label{eq:Ijl}
I(\theta_s)=I(\varepsilon_{s-1})=j_{s-1}.
\end{equation}

We use the following results on $j_l$.
\begin{itemize}
\item \citet{Ennola} conjectured that $j_l=1$ for every $l\ge3$, proved
it for $3\le l\le500$, and showed that $j_l$ is coprime to $30$.
\item \citet{Lou17} showed that $j_l$ is coprime to $19!$ for every
$l\ge3$, and that $j_l=1$ for $3\le l\le5\cdot10^7$.
\item \citet{ChoiKim}, proving conjectures of \citet{Lou21}, showed that
for every prime $p\ge3$ only finitely many $l\ge3$ have $p\mid j_l$.
\end{itemize}
Ennola's conjecture itself remains open.

\begin{theorem}[Coefficients in the second family]\label{thm:index}
Let $K$ be a totally real non-Galois cubic field, and let $z_0,z_1$ be
exceptional units of $K$ with $z_0,1-z_0$ independent and
\[
m\,\delta(z_0)+\delta(z_1)=0,\qquad m\in\Z_{\ge2}.
\]
Then there are parameters $l_0,l_1\ge3$ with
$K=\Q(\varepsilon_{l_0})=\Q(\varepsilon_{l_1})$, $z_i$ anharmonically
equivalent to $\varepsilon_{l_i}$, and
\[
j_{l_1}=m\,j_{l_0}.
\]
Consequently:
\begin{enumerate}
\item $l_1>5\cdot10^7$;
\item every prime divisor of $m$ is at least $23$; in particular
$m\ge23$;
\item for each $m$, at most finitely many fields $K$ carry such a
relation;
\item if Ennola's conjecture holds, no such relation exists.
\end{enumerate}
\end{theorem}

\begin{proof}
By Lemma~\ref{lem:generatedorder}, each $z_i$ has an anharmonic transform
whose minimal polynomial is $f_t$ with $t\ge-1$ or $g_s$ with $s\ge3$.
The fields of the first family are cyclic, and $g_3$ gives the
conductor-$7$ field, so $K$ non-Galois forces $g_s$ with $s\ge4$.  By
$-g_s(1-X)=F_{s-1}(X)$, $z_i$ is anharmonically equivalent to
$\varepsilon_{l_i}$ with $l_i=s_i-1\ge3$, and
$K=\Q(z_i)=\Q(\varepsilon_{l_i})$.  Lemma~\ref{lem:indexlaw} and
\eqref{eq:Ijl} give $j_{l_1}=I(z_1)=m\,I(z_0)=m\,j_{l_0}$.  Hence
$j_{l_1}\ge2$, and \citet{Lou17} gives (1).  Every prime divisor of $m$
divides $j_{l_1}$, which is coprime to $19!$; this is (2).  For (3), fix
a prime $p\mid m$; by \citet{ChoiKim} only finitely many $l_1$ have
$p\mid j_{l_1}$, and $l_1$ determines $K$.  Part (4) is immediate.
\end{proof}

Theorem~\ref{thm:secondempty} sharpens Theorem~\ref{thm:index} for
arbitrary pairs of exceptional units.

\begin{corollary}[Both indices]\label{cor:index}
In Theorem~\ref{thm:index} also $j_{l_0}\ge23$. Consequently both
parameters exceed $5\cdot10^7$, and $j_{l_1}=m\,j_{l_0}\ge23^2$.
\end{corollary}

\begin{proof}
Suppose $j_{l_0}=1$, and let $s_0=l_0+1$, so that $z_0$ is anharmonically
equivalent to $\theta_{s_0}$. Then $\mathcal O_K^\times=\langle-1,
\theta_{s_0},1-\theta_{s_0}\rangle\subset\Z[\theta_{s_0}]$, so $z_1$ and
$1-z_1$ are units of $\Z[\theta_{s_0}]$, and by
Theorem~\ref{thm:secondempty} $z_1$ lies in the anharmonic orbit of
$\theta_{s_0}$. Then $I(z_1)=I(z_0)$, contradicting $I(z_1)=m\,I(z_0)$
with $m\ge2$. So $j_{l_0}\ge2$; since $j_{l_0}$ is coprime to $19!$
\citep{Lou17}, $j_{l_0}\ge23$, and $j_{l_0}=1$ for $l_0\le5\cdot10^7$
gives $l_0>5\cdot10^7$.
\end{proof}

\begin{remark}\label{rem:indexscope}
Since Theorem~\ref{thm:index} and Corollary~\ref{cor:index} only use
exceptional units and the wedge relation, they apply verbatim to the
four-unit points of rank-two Nahm systems with symmetrizer
$\operatorname{diag}(1,m)$, for which the Nahm equations give
$m\,\delta(x)+\delta(y)=0$.
\end{remark}

\section[Complex cubic fields]{Cubic fields of signature \texorpdfstring{$(1,1)$}{(1,1)}}\label{sec:complex}

Here the unit rank is one.  Hence every exceptional unit satisfies
\[
\delta(z)=0,
\]
so the determinant law no longer determines the coefficient in
\eqref{eq:twoterm}.  We therefore treat this rank-one branch directly in
the five-term/reflection calculus of Definition~\ref{def:certified}.

The signature-$(1,1)$ clause of Theorem~\ref{thm:nagell} leaves exactly the
fields of discriminants $-23$ and $-31$.  Their exceptional units can be
enumerated exactly.

\begin{remark}[The endpoint unit test]\label{rem:complexendpoint}
There is a simple preliminary obstruction behind this very short list.
If $\alpha$ is an exceptional algebraic integer of degree three with monic
minimal polynomial $q(X)$, then
\[
|q(0)|=|N_{K/\Q}(\alpha)|=1,
\qquad
|q(1)|=|N_{K/\Q}(1-\alpha)|=1.
\]
Thus most cubic generators are excluded before any Rogers-dilogarithm
calculation is needed.

For example, let $\alpha\in(0,1)$ be the real root of
\[
q(X)=X^3+2X-1.
\]
Then
\[
q(0)=-1,\qquad q(1)=2,
\]
and hence
\[
|N_{K/\Q}(1-\alpha)|=2.
\]
Thus $1-\alpha$ is not a unit, so $\alpha$ is not an exceptional unit.
Moreover
\[
\operatorname{disc}(q)=-59.
\]
Since $-59$ is squarefree, this is the field discriminant; consequently
Theorem~\ref{thm:nagell} shows more strongly that this cubic field contains
no exceptional unit at all.  The same endpoint test $|q(0)|=|q(1)|=1$
is the elementary norm condition underlying the two exceptional cubic
families used throughout the paper.
\end{remark}

\begin{proposition}[Exceptional units in the two complex cubic fields]
\label{prop:excunits}
Let $\eta$ denote the fundamental unit normalized to lie in $(0,1)$ at the
real place.  Then the exceptional units in $(0,1)$ are
\[
\eta,\eta^2,\eta^3,\eta^5
\quad(\operatorname{disc}K=-23),
\qquad
\eta,\eta^3
\quad(\operatorname{disc}K=-31).
\]
They satisfy
\[
\eta+\eta^5=1,\qquad \eta^2+\eta^3=1
\]
in discriminant $-23$, and
\[
\eta+\eta^3=1
\]
in discriminant $-31$.
\end{proposition}

\begin{proof}
For discriminant $-23$, $\eta$ is the real root in $(0,1)$ of
\[
X^3+X^2-1=0;
\]
for discriminant $-31$, it is the root of
\[
X^3+X-1=0.
\]
We first justify that these roots are fundamental units. Both roots exceed
$2/3$, so $\alpha=\eta^{-1}$ satisfies $1<\alpha<3/2<25/16$.
If $\alpha=\beta^k$ for a positive unit $\beta>1$ and an integer
$k\ge2$, then $1<\beta<5/4$. The norm of $\beta$ is $1$, and its
two nonreal conjugates have modulus $\beta^{-1/2}$. Write its minimal
polynomial as $X^3-AX^2+BX-1$, with $A,B\in\Z$. The trace bound gives
$-1<A<4$, hence $A\in\{0,1,2,3\}$, while
\[
B=A\beta-\beta^2+\beta^{-1}.
\]
For each of these four values of $A$, the function
$h_A(t)=At-t^2+t^{-1}$ is strictly decreasing on $1<t\le5/4$,
since $2t+t^{-2}>3\ge A$. Moreover
\[
h_A(1)=A,\qquad h_A(5/4)=\frac{5A}{4}-\frac{61}{80}>A-1.
\]
Thus $A-1<B<A$, contradicting $B\in\Z$. It follows that $\alpha$,
and hence $\eta$, generates the positive unit group.

Since the unit rank is one, every positive unit at the real place is a
power of $\eta$.  If $\xi,\overline\xi$ are the complex conjugates, then
$|\xi|=\eta^{-1/2}$ and
\[
|N(1-\eta^k)|
=(1-\eta^k)|1-\xi^k|^2
\ge
(1-\eta^k)(\eta^{-k/2}-1)^2.
\]
Writing $t=\eta^{-k/2}>1$, the right side is
$(t-1)^3(t+1)/t^2$, which is increasing in $t$.  The exact isolating
bounds
\[
0.754<\eta<0.755\quad(\operatorname{disc}=-23),
\qquad
0.682<\eta<0.683\quad(\operatorname{disc}=-31)
\]
show that the lower bound exceeds $1$ for $k\ge6$ and $k\ge4$,
respectively.  Exact norms for the remaining exponents are
\[
1,1,1,5,1
\quad\text{and}\quad
1,3,1,
\]
giving the stated lists.  The displayed complement identities follow from
the two cubic equations.
\end{proof}

The following classical observation rules out a hidden one-term collapse in
both fields.

\begin{proposition}[Single-symbol obstruction, after Zagier]\label{prop:singlesymbol}
Let $x\in\overline{\Q}\cap(0,1)$. If $[x]$ represents a torsion
element of $\Bx(\mathbb C)$, equivalently if its image vanishes in
$\Bx(\mathbb C)\otimes\Q$, then
\[
x\in\left\{\frac12,\ \frac{\sqrt5-1}{2},\ \frac{3-\sqrt5}{2}\right\}.
\]
In particular, the fundamental units occurring in Proposition~\ref{prop:excunits}
have nonzero classes in $\Bx(\mathbb C)\otimes\Q$.
\end{proposition}

\begin{proof}
This is the real $(0,1)$ part of the single-symbol torsion classification
in $\Bx(\mathbb C)$ discussed by \citet[\S2B(c)]{Zagier07}: up to the six anharmonic
transformations, the only torsion symbols of the form $[x]$ are represented
by $x=1/2$ and the inverse golden ratio.  Their representatives in $(0,1)$
are exactly the three values displayed above.  The two cubic units in
Proposition~\ref{prop:excunits} have degree three and therefore are not among
them.
\end{proof}

Lewin's identity itself follows from one specialization of the five-term
identity and one reflection.  Setting $x=y=\eta$ in
\eqref{eq:fiveterm} gives Landen's duplication formula
\begin{equation}\label{eq:lewinlanden}
2\Lr(\eta)=\Lr(\eta^2)+2\Lr\!\left(\frac{\eta}{1+\eta}\right).
\end{equation}
For discriminant $-23$ one has $\eta^2(1+\eta)=1$, hence
$\eta/(1+\eta)=\eta^3$, while $\eta^2+\eta^3=1$.  Therefore
\[
2\Lr(\eta)=\Lr(\eta^2)+2\Lr(\eta^3)
=2\Lr(1)-\Lr(\eta^2),
\]
and so
\begin{equation}\label{eq:plastic}
2\Lr(\eta)+\Lr(\eta^2)=2\Lr(1)=\frac{\pi^2}{3}.
\end{equation}
This is Lewin's plastic identity \citep{Lewin82}.

\begin{theorem}[Complex cubic certified classification]\label{thm:lewin}
Let $K$ be a cubic field of signature $(1,1)$ and let
$z_0,z_1\in K\cap(0,1)$ be exceptional units.  Among five-term certified
relations
\[
m\,\Lr(z_0)+\Lr(z_1)=c\pi^2,
\qquad m\in\Z_{>1},\quad c\in\Q,
\]
the only one, up to simultaneous reflection, is Lewin's identity
\eqref{eq:plastic} in discriminant $-23$.  The field of discriminant $-31$
contributes only the $m=1$ reflection identity.
\end{theorem}

\begin{proof}
Work first in $\mathcal P(\mathbb C)\otimes\Q$, so that auxiliary
certificate arguments outside $K$ are also included. The five-term
relations hold there, and reflection gives
\[
[x]+[1-x]=0.
\]
Any five-term certified two-term Rogers relation therefore gives
\begin{equation}\label{eq:certbloch}
m[z_0]+[z_1]=0.
\end{equation}
Since $z_0,z_1$ are exceptional units in a field of unit rank one,
$\delta(z_0)=\delta(z_1)=0$. Thus their symbols, and the resulting
relation, lie in $\Bx(\mathbb C)\otimes\Q$.

For discriminant $-31$, Proposition~\ref{prop:excunits} gives only
$\eta$ and $\eta^3=1-\eta$, so
\[
[\eta^3]=-[\eta].
\]
If \eqref{eq:certbloch} held with $m>1$, it would force a non-zero integer
multiple of $[\eta]$ to vanish.  Proposition~\ref{prop:singlesymbol} rules
this out.  Thus only $m=1$, namely reflection, occurs.

For discriminant $-23$, reflection gives
\[
[\eta^5]=-[\eta],\qquad [\eta^3]=-[\eta^2].
\]
The five-term specialization underlying \eqref{eq:lewinlanden} gives
\[
2[\eta]=[\eta^2]+2[\eta^3],
\]
so
\[
[\eta^2]=-2[\eta],\qquad
[\eta^3]=2[\eta],\qquad
[\eta^5]=-[\eta].
\]
By Proposition~\ref{prop:singlesymbol}, $[\eta]\ne0$ after tensoring with
$\Q$.  Hence the four exceptional units have coefficients
\[
1,\ -2,\ 2,\ -1
\]
relative to the non-zero class $[\eta]$.  Inspecting
$m[z_0]+[z_1]=0$, the only positive integer ratio $m>1$ is $m=2$, arising
from $(z_0,z_1)=(\eta,\eta^2)$ and its simultaneous reflection
$(\eta^5,\eta^3)$.  The corresponding Rogers evaluation is
\eqref{eq:plastic}.
\end{proof}

\begin{remark}[Exact scope]
Theorem~\ref{thm:lewin} is a classification in the five-term certified class
of Definition~\ref{def:certified}.  It makes no assertion about a hypothetical
numerical equality which has no five-term/reflection certificate.  No
period non-vanishing statement enters the proof.
\end{remark}

\appendix

\section{Proof of finite-prime rigidity}\label{app:forcing}
\subsection{The finite-valuation fan and its symmetry}

Put $y=1-x$, $V=1-U$, and let $\mathfrak p$ be a finite prime of $K$.
We normalize $v_{\mathfrak p}$ to have value group $\Z$.

\begin{lemma}[Valuation fan]\label{lem:fan}
For every finite prime $\mathfrak p$,
\[
\bigl(\vp(x),\vp(y)\bigr)
\in
\R_{>0}(1,0)\ \cup\ \R_{>0}(0,1)\ \cup\
\R_{>0}(-1,-1)\ \cup\ \{(0,0)\}.
\]
The same assertion holds for $(\vp(U),\vp(V))$.
\end{lemma}

\begin{proof}
If $\vp(x)>0$, then $\vp(1-x)=0$.  If $\vp(x)<0$, the ultrametric
inequality gives $\vp(1-x)=\vp(x)$.  Finally, if $\vp(x)=0$, then
$\vp(1-x)\ge0$.  These are precisely the three stated rays and the
origin.  The proof for $U+V=1$ is identical.
\end{proof}

Write
\[
R_x=\R_{>0}(1,0),\qquad
R_y=\R_{>0}(0,1),\qquad
R_\infty=\R_{>0}(-1,-1).
\]
The exponent matrix satisfies
\begin{equation}\label{eq:valuationmatrix}
\begin{pmatrix}\vp(U)\\ \vp(V)\end{pmatrix}
=M\begin{pmatrix}\vp(x)\\ \vp(y)\end{pmatrix}.
\end{equation}
Thus every source ray that occurs is sent to one of the same three
target rays.

The six standard transformations of $X+Y=1$ are
\begin{equation}\label{eq:s3pairs}
(X,Y),\quad (Y,X),\quad
(X^{-1},-Y/X),\quad (-Y/X,X^{-1}),\quad
(Y^{-1},-X/Y),\quad (-X/Y,Y^{-1}).
\end{equation}
Each transformed coordinate is a signed monomial in $X$ and $Y$.
The signs are invisible to finite valuations, while the exponent
vectors are acted on by unimodular matrices.  These transformations
permute $R_x,R_y,R_\infty$ as the full symmetric group $S_3$.

\begin{lemma}[Two-ray normalization]\label{lem:tworaynormal}
Assume that $\det M\ne0$ and that at least two distinct source rays
occur for $(x,y)$.  Then independent transformations
\eqref{eq:s3pairs} of the source solution and the target solution
reduce the exponent matrix to
\[
\begin{pmatrix}A&0\\0&B\end{pmatrix},
\qquad A,B\in\Z_{>0},\qquad AB=|\det M|.
\]
Consequently a rational M\"obius transform $X$ of $x$ satisfies
\begin{equation}\label{eq:tworayeq}
\epsilon X^A+\eta(1-X)^B=1
\end{equation}
for some $\epsilon,\eta\in\{\pm1\}$.
\end{lemma}

\begin{proof}
Choose two occurring source rays.  Since $M$ is nonsingular, their
images are distinct target rays.  Let $S,T\in\mathrm{GL}_2(\Z)$ be,
respectively, the matrices of a source $S_3$ transformation and an
independent target $S_3$ transformation that send these two pairs of
rays to the coordinate rays in the same order.  If
$\mathbf v'=S\mathbf v$ and $\mathbf w'=T\mathbf w$, then
\[
\mathbf w'=TMS^{-1}\mathbf v',
\]
so the transformed exponent matrix is
\[
M'=TMS^{-1}.
\]
Because $S$ and $T$ are unimodular,
\[
|\det M'|=|\det M|=|D|.
\]
The two columns of $M'$ lie on the positive coordinate rays, hence
\[
M'=\begin{pmatrix}A&0\\0&B\end{pmatrix},
\qquad A,B\in\Z_{>0}.
\]
It follows immediately that
\[
AB=|\det M'|=|D|.
\]
The transformed source and target pairs still solve $X+Y=1$.
Because the maps in \eqref{eq:s3pairs} are monomial only up to sign,
the signs in \eqref{eq:tworayeq} record the possible factors $-1$;
they do not affect the valuation matrices or determinant calculation.
\end{proof}

\subsection{Five real-root lemmas}

The following elementary lemmas are the archimedean input.  Their
statements are arranged so that possible rational boundary roots such
as $-1$, $0$, or $1$ do not affect the exclusion of an irreducible
totally real cubic factor.

\begin{lemma}[Two-ray polynomial]\label{lem:rootA}
Let $A,B\in\Z_{>0}$ with $AB>1$, and let
$\epsilon,\eta\in\{\pm1\}$.  The polynomial
\[
f(X)=\epsilon X^A+\eta(1-X)^B-1
\]
has no root in $(0,1)$ and at most one real root outside
$\{0,1\}$.  In particular, it has no irreducible totally real cubic
factor.
\end{lemma}

\begin{proof}
For $0<X<1$, one has $X^A\le X$ and $(1-X)^B\le1-X$.  In the
$+,+$ case equality in both inequalities would force $A=B=1$, excluded
by $AB>1$; the other sign cases are immediate.  Hence $f<0$ on
$(0,1)$.

Let $X=1+t$ with $t>0$ and put $\sigma=\eta(-1)^B$.  The equal-sign
cases
\[
\epsilon(1+t)^A+\sigma t^B-1
\]
never vanish.  In the mixed cases one obtains either
\[
G_{A,B}(t)=(1+t)^A-t^B-1
\]
or
\[
H_{A,B}(t)=t^B-(1+t)^A-1.
\]
If $B\le A$, then $G_{A,B}(t)>0$ and $H_{A,B}(t)<0$ for $t>0$.
If $B>A$, the nonzero coefficient sequence of each polynomial has
exactly one sign change, while its signs at $0^+$ and $+\infty$ are
opposite.  Each therefore has exactly one positive root.  Thus a root
with $X>1$ is possible only when $B>A$, and it is unique.

The substitution $X\mapsto1-X$ interchanges $A$ and $B$ and swaps the
two exterior intervals.  Hence a root with $X<0$ is possible only when
$A>B$, and it too is unique.  The two possibilities are mutually
exclusive.
\end{proof}

For $s\ge1$, put
\[
S_s(Y)=1+Y+\cdots+Y^{s-1}.
\]

\begin{lemma}[First single-ray polynomial]\label{lem:rootB}
Let $A\in\Z_{\ge0}$ and $Q,s\in\Z_{>0}$, with $(A+1)s>1$.  The
polynomial
\[
F(Y)=(1-Y)^A-Y^Q S_s(Y)
\]
has exactly one root in $(0,1)$, no root in $(-1,0)$, and at most one
root in $(-\infty,-1)\cup(1,\infty)$.  Hence it has at most two real
roots other than the possible rational root $Y=-1$.
\end{lemma}

\begin{proof}
On $(0,1)$, the first term is nonincreasing and the second is strictly
increasing; their endpoint values give exactly one intersection.  If
$Y=-z$ with $0<z<1$, then $S_s(-z)>0$ and $S_s(-z)\le1$.  If $Q$ is
odd, both terms of $F(-z)$ are positive; if $Q$ is even, then
\[
z^Q S_s(-z)<1\le(1+z)^A.
\]
Thus there is no root in $(-1,0)$.

Set $N=Q+s-1$.  For $Y>1$, a root is possible only when $A$ is even.
If $A=0$, then
\[
Y^Q S_s(Y)>1=(Y-1)^A,
\]
so no root exists.  Assume for the remainder of this positive-axis
case that $A\ge1$.  If $A\le N$, then
\[
(Y-1)^A<Y^A\le Y^N\le Y^Q S_s(Y),
\]
so no root exists.  If $A>N$, the quotient
\[
R_+(Y)=\frac{(Y-1)^A}{Y^Q S_s(Y)}
\]
is strictly increasing, since
\[
Y\frac{R_+'(Y)}{R_+(Y)}
=\frac{AY}{Y-1}-Q-Y\frac{S_s'(Y)}{S_s(Y)}
>A-Q-(s-1)=A-N>0.
\]
It therefore crosses $1$ at most once.

For $Y=-z<-1$, a root is possible only in the parity-compatible case,
when
\[
(1+z)^A=z^Q|S_s(-z)|.
\]
If $A\ge N$, then
\[
z^Q|S_s(-z)|\le z^N<(1+z)^A,
\]
so no root exists.  If $A<N$, put
\[
R_-(z)=\frac{z^Q|S_s(-z)|}{(1+z)^A}.
\]
For odd $s$ one has
$|S_s(-z)|=(1+z^s)/(1+z)$, and
\[
z\frac{R_-'(z)}{R_-(z)}
=Q+\frac{s z^s}{1+z^s}-\frac{(A+1)z}{1+z}.
\]
Put $u=z/(1+z)\in(0,1)$.  Since
\[
\frac{z^s}{1+z^s}\ge u,
\]
we obtain
\[
z\frac{R_-'(z)}{R_-(z)}
\ge Q+(s-A-1)u.
\]
If $s-A-1\ge0$, this is at least $Q>0$.  If $s-A-1<0$, then
$u<1$ gives
\[
Q+(s-A-1)u>Q+s-A-1=N-A>0.
\]
Thus $R_-'(z)>0$ in either case.  For even $s$ one has
$|S_s(-z)|=(z^s-1)/(1+z)$, and the same conclusion follows from
\[
Q+\frac{s z^s}{z^s-1}-\frac{(A+1)z}{1+z}>N-A>0.
\]
Thus there is at most one root below $-1$.  The two exterior cases are
mutually exclusive, since they require $A>N$ and $A<N$, respectively.
\end{proof}

\begin{lemma}[Second single-ray polynomial]\label{lem:rootC}
Let $1\le q<s$ and $H\ge1$.  The polynomial
\[
G(Y)=Y^q+Y^s-(1-Y)^H
\]
has exactly one root in $(0,1)$, no root in $(-1,0)$, and at most one
root in $(-\infty,-1)\cup(1,\infty)$.  Hence it has at most two real
roots other than the possible rational root $Y=-1$.
\end{lemma}

\begin{proof}
The assertion on $(0,1)$ follows because $Y^q+Y^s$ is strictly
increasing and $(1-Y)^H$ is strictly decreasing.  For $Y=-z$ with
$0<z<1$,
\[
|(-z)^q+(-z)^s|\le z^q+z^s<1+z\le(1+z)^H,
\]
so no root occurs.

For $Y>1$, a root requires $H$ even.  If $H\le s$, then
\[
(Y-1)^H<Y^H\le Y^s<Y^q+Y^s.
\]
If $H>s$, the quotient
\[
\frac{(Y-1)^H}{Y^q+Y^s}
\]
is strictly increasing, because its logarithmic derivative is greater
than $H/(Y-1)-s/Y>0$.  Hence there is at most one root above $1$.

For $Y=-z<-1$, positivity requires $s$ even.  The positive left side is
then either $z^s+z^q$ or $z^s-z^q$.  If $H\ge s$, the binomial
expansion gives
\[
(1+z)^H>z^s+z^q>z^s-z^q,
\]
so no root exists.  Suppose $H<s$ and put $d=s-q$.  For the plus sign, put $u=z/(1+z)\in(0,1)$.  Since
$z^d/(1+z^d)\ge u$,
\[
z\frac{d}{dz}\log\frac{z^q(1+z^d)}{(1+z)^H}
=q+\frac{d z^d}{1+z^d}-\frac{Hz}{1+z}
\ge q+(d-H)u.
\]
If $d-H\ge0$, this is at least $q>0$.  If $d-H<0$, then $u<1$ gives
\[
q+(d-H)u>q+d-H=s-H>0.
\]
For the minus sign the same logarithmic derivative, with
$1+z^d$ replaced by $z^d-1$, is even larger because
$z^d/(z^d-1)>1$.
Thus there is at most one root below $-1$.  The two exterior cases
require $H>s$ and $H<s$, respectively, and are mutually exclusive.
\end{proof}

\begin{lemma}[Denominator-to-numerator polynomial]\label{lem:rootD}
Let $A,B,H\in\Z_{>0}$.  The polynomial
\[
P(T)=T^A(1+T)^H+T^B-1
\]
has at most two real roots other than the possible rational root
$T=-1$.  In particular, it has no irreducible totally real cubic
factor.
\end{lemma}

\begin{proof}
The left side of
\[
T^A(1+T)^H+T^B=1
\]
is strictly increasing from $0$ to $+\infty$ on $(0,\infty)$, so there
is exactly one positive root.  If $T=-z$ with $0<z<1$, then
\[
\bigl|(-z)^A(1-z)^H+(-z)^B\bigr|
\le z(1-z)+z=2z-z^2<1,
\]
so there is no root in $(-1,0)$.

Let $T=-z<-1$ and put $C=A+H$.  Then
\[
P(-z)=(-1)^C z^A(z-1)^H+(-1)^Bz^B-1.
\]
If $C$ and $B$ are both even, the two positive terms cannot sum to
$1$, since $z^B>1$.  If $C$ and $B$ are both odd, the left-hand side
is negative.  The only remaining, parity-compatible equations are
\[
z^A(z-1)^H=1+z^B
\]
and
\[
z^B-1=z^A(z-1)^H.
\]
For the first, no root exists when $C\le B$; when $C>B$, the quotient
$z^A(z-1)^H/(1+z^B)$ is strictly increasing because its logarithmic
derivative is greater than $(C-B)/z>0$.

For the second, divide by $z-1$:
\[
1+z+\cdots+z^{B-1}=z^A(z-1)^{H-1}.
\]
No root exists when $C\le B$.  When $C>B$ and $B=1$, the left side is
$1$, while the logarithmic derivative of the right side is
\[
\frac{A}{z}+\frac{H-1}{z-1}>0.
\]
When $C>B$ and $B\ge2$, the quotient of the right side by the left
side is strictly increasing, since
\[
\frac{z(1+2z+\cdots+(B-1)z^{B-2})}
     {1+z+\cdots+z^{B-1}}<B-1.
\]
Thus at most one root lies below $-1$.
\end{proof}

\begin{lemma}[Residual denominator polynomial]\label{lem:rootE}
Let $A,H,K\in\Z_{>0}$.  The polynomial
\[
Q(T)=1+T^K-T^A(1+T)^H
\]
has at most two real roots other than the possible rational root
$T=-1$.  In particular, it has no irreducible totally real cubic
factor.
\end{lemma}

\begin{proof}
There is no root in $(-1,0)$.  Indeed, write $T=-z$ with $0<z<1$.
All parity cases are immediate except $A$ even and $K$ odd, where a
root would imply
\[
1+z+\cdots+z^{K-1}=z^A(1-z)^{H-1};
\]
the left side is at least $1$ and the right side is less than $1$.

Put $C=A+H$.  On the positive axis, the coefficient sequence of $Q$
has at most one sign change when $K\le C$ and at most two when $K>C$.
Descartes' rule therefore gives at most one and at most two positive
roots, respectively.

For $T=-z<-1$, a root can occur only when $K$ and $C$ have the same
parity.  The resulting equations are
\[
1+z^K=z^A(z-1)^H
\]
or
\[
1+z+\cdots+z^{K-1}=z^A(z-1)^{H-1}.
\]
If $K>C$, neither can hold, by comparison with $z^K$ and $z^{K-1}$.
If $K\le C$, each has at most one solution.  For the first, the
logarithmic derivative of $z^A(z-1)^H/(1+z^K)$ is positive.  For the
second, if $K=1$ the left side is $1$ and the logarithmic derivative
of the right side is
\[
\frac{A}{z}+\frac{H-1}{z-1}>0.
\]
If $K\ge2$, the quotient of the right side by the left side is strictly
increasing because
\[
\frac{z(1+2z+\cdots+(K-1)z^{K-2})}
     {1+z+\cdots+z^{K-1}}<K-1.
\]
Thus, if $K\le C$, there is at most one positive root and at most one
root below $-1$; if $K>C$, there are at most two positive roots and no
root below $-1$.
\end{proof}

\subsection{Completion of the proof}

\paragraph{Proof architecture.}
By Lemma~\ref{lem:fan}, every nonzero finite valuation vector of
$(x,1-x)$ lies on one of the three rays $R_x,R_y,R_\infty$.
Lemma~\ref{lem:tworays} shows that at most one of these rays can occur.
Lemma~\ref{lem:numerator} excludes $R_x$ and $R_y$, while
Lemma~\ref{lem:denominator} excludes $R_\infty$.  The finite support is
therefore empty.

We now exclude every possible nonempty finite support.

\begin{lemma}[Two support rays are impossible]\label{lem:tworays}
At most one of $R_x,R_y,R_\infty$ occurs among the nonzero valuation
vectors $(\vp(x),\vp(y))$.
\end{lemma}

\begin{proof}
If two rays occur, Lemma~\ref{lem:tworaynormal} gives a rational
M\"obius transform $X$ of $x$ satisfying \eqref{eq:tworayeq} with
$AB=|D|>1$.  By Lemma~\ref{lem:mobius}, the element $X$ has degree
three and three distinct real nonrational conjugates, all of which
satisfy \eqref{eq:tworayeq}.  This contradicts
Lemma~\ref{lem:rootA}.
\end{proof}

\begin{lemma}[A numerator ray is impossible]\label{lem:numerator}
The sole occurring finite-valuation ray cannot be $R_x$ or $R_y$.
\end{lemma}

\begin{proof}
By symmetry, suppose that the sole ray is $R_x$.  At a supporting
prime, the first column $(p,r)^T$ of $M$ lies on one of the target rays.
It is nonzero because $\det M\ne0$.  After interchanging $U$ and $V$,
there are two cases.

First suppose $p>0$ and $r=0$.  Then
\[
U=x^py^q,\qquad V=y^s.
\]
Since $0<y,V<1$, one has $s>0$, and $|D|=ps>1$.  If $q\ge0$, then
$U\le x$ and $V\le y$, with at least one strict inequality, contrary
to $U+V=x+y=1$.  Thus $q=-Q$ with $Q>0$.  Starting from
\[
x^py^{-Q}+y^s=1,
\]
multiply by $y^Q$ to obtain
\[
x^p+y^{Q+s}=y^Q.
\]
Using $x=1-y$ gives
\[
(1-y)^p=y^Q(1-y^s).
\]
Since
\[
1-y^s=(1-y)S_s(y),
\]
division by $1-y\ne0$ yields
\[
(1-y)^{p-1}=y^Q S_s(y).
\]
Lemma~\ref{lem:rootB} shows that this polynomial has at most two
nonrational real roots, whereas Lemma~\ref{lem:mobius}, applied to
$y=1-x$, shows that the three conjugates of the totally real cubic
element $y$ are distinct and nonrational.  This is a contradiction.

The remaining target ray is $R_\infty$.  Then $p=r=-H<0$.  Moreover
$q,s>0$: if, for example, $q\le0$, then $x^{-H}>1$ and $y^q\ge1$, so
$U=x^{-H}y^q>1$, contrary to $0<U<1$; the argument for $s$ is the same.
Since $D=H(q-s)\ne0$, interchange the summands so that $q<s$.  Multiplication
by $x^H$ gives
\[
y^q+y^s=(1-y)^H.
\]
By Lemma~\ref{lem:mobius}, $y=1-x$ has three distinct real
nonrational conjugates, contradicting Lemma~\ref{lem:rootC}.
\end{proof}

\begin{lemma}[The common-denominator ray is impossible]
\label{lem:denominator}
The sole occurring finite-valuation ray cannot be $R_\infty$.
\end{lemma}

\begin{proof}
If only $R_\infty$ occurs, put
\[
t=\frac{x}{y}.
\]
At every finite prime in the support one has
$\vp(x)=\vp(y)$, while outside the support both valuations vanish.
Hence
\[
\vp(t)=\vp(x)-\vp(y)=0
\]
for every finite prime $\mathfrak p$, and therefore
\[
t\in\OK^\times.
\]
This fact is not needed below; it only records that the
common-denominator case is already unit-like in the quotient.  Moreover
\[
x=\frac{t}{1+t},\qquad y=\frac1{1+t}.
\]
By Lemma~\ref{lem:mobius}, the rational M\"obius transform $t=x/y$
has degree three and three distinct real conjugates.
In particular $t\ne1$, since $t=1$ would give $x=y=1/2$.

Suppose first that $R_\infty$ maps to a numerator ray.  After
interchanging $U,V$,
\[
p+q=-H<0,\qquad r+s=0.
\]
The equation $U+V=1$ is
\[
t^p(1+t)^H+t^r=1.
\]
If $0<t<1$, then each summand lies in $(0,1)$, which forces
$p,r>0$, so this is already of the form
\[
t^A(1+t)^H+t^B=1
\]
with $A=p$ and $B=r$.  If $t>1$, positivity instead forces $p,r<0$.
Put $T=t^{-1}\in(0,1)$.  Then
\[
t^p(1+t)^H=T^{-p-H}(1+T)^H,
\qquad t^r=T^{-r}.
\]
Because the first summand is $<1$, one has $-p-H>0$; also $-r>0$.
Thus the transformed equation is
\[
T^A(1+T)^H+T^B=1,
\qquad A=-p-H>0,
\qquad B=-r>0.
\]
In either case Lemma~\ref{lem:rootD} is contradicted.

It remains that $R_\infty$ maps to itself.  Then
\[
p+q=r+s=H>0
\]
and
\begin{equation}\label{eq:residual}
t^p+t^r=(1+t)^H.
\end{equation}
Moreover
\[
D=H(p-r),
\]
so $p\ne r$.  Assume $p>r$.  If $r<0$, multiplication of
\eqref{eq:residual} by $t^{-r}$ gives
\[
1+t^{p-r}=t^{-r}(1+t)^H.
\]
Thus this case has the form
\[
1+t^K=t^A(1+t)^H,
\qquad K=p-r>0,
\qquad A=-r>0.
\]

Suppose now that $r\ge0$.  If $p\le H$, the binomial expansion of
$(1+t)^H$ strictly dominates $t^r+t^p$ for every $t>0$, except for
$H=1$ and $(r,p)=(0,1)$, which gives $|D|=1$.  Hence a nontrivial
configuration must have $p>H$.  Put $T=t^{-1}$.  Dividing
\eqref{eq:residual} by $t^p$ gives
\[
1+T^{p-r}=T^{p-H}(1+T)^H.
\]
Thus
\[
K=p-r>0,
\qquad A=p-H>0.
\]
In either case Lemma~\ref{lem:rootE} is contradicted (with $T$ renamed
as $t$ in the second case).
\end{proof}

\begin{proof}[Proof of Theorem~\ref{thm:forcing}]
If the finite valuation support of $(x,1-x)$ were nonempty,
Lemma~\ref{lem:tworays} would leave exactly one support ray.
Lemma~\ref{lem:numerator} excludes $R_x$ and $R_y$, and
Lemma~\ref{lem:denominator} excludes $R_\infty$.  Hence the support is
empty, so
\[
\vp(x)=\vp(1-x)=0
\]
for every finite prime $\mathfrak p$.  Therefore
\[
x,\ 1-x\in\OK^\times.
\]
Finally,
\[
U=x^p(1-x)^q,
\qquad
1-U=x^r(1-x)^s
\]
with integral exponents.  Both are products of integral powers of
units, so negative exponents cause no difficulty and
\[
U,\ 1-U\in\OK^\times.
\]
\end{proof}

\section{Details of the effective bound for the second family}\label{app:effective}
\subsection{Explicit root intervals}
The estimates below are used with explicit constants, so we record
them as inequalities rather than as asymptotics.

\begin{lemma}[Root intervals]\label{lem:rootint}
For $s\ge10$ the roots of $g_s$ satisfy
\[
s-1<\theta_1<s,\qquad
1+\tfrac1s<\theta_2<1+\tfrac2s,\qquad
-\tfrac1s<\theta_3<-\tfrac1{s^{2}} .
\]
Consequently, with $\tau_i=1-\theta_i$,
\[
\bigl|\log|\theta_1|-\log s\bigr|<\tfrac1{s-1},\quad
\bigl|\log|\tau_1|-\log s\bigr|<\tfrac2{s-2},\quad
\tfrac1{s+1}<\log|\theta_2|<\tfrac2s,
\]
\[
\log s-\log2<\bigl|\log|\tau_2|\bigr|<\log s,\quad
\log s<\bigl|\log|\theta_3|\bigr|<2\log s,\quad
0<\log|\tau_3|<\tfrac2s .
\]
\end{lemma}

\begin{proof}
Each root interval follows from evaluating $g_s$ at its endpoints and
applying the intermediate value theorem; for instance
$g_s(s-1)=-s^{2}+3s-1<0<1=g_s(s)$ for $s\ge3$. For the second root
the sign change runs the other way:
\[
g_s\Bigl(1+\tfrac1s\Bigr)=\frac{2s+1}{s^{3}}>0,
\qquad
g_s\Bigl(1+\tfrac2s\Bigr)=-\frac{s^{3}+2s^{2}-8s-8}{s^{3}}<0
\quad(s\ge10),
\]
whence $1+\tfrac1s<\theta_2<1+\tfrac2s$. For the logarithmic
statements: $s-1<\theta_1<s$ gives
$\log s+\log(1-\tfrac1s)<\log\theta_1<\log s$, and
$\log(1-\tfrac1s)>-\tfrac1{s-1}$; likewise $|\tau_1|=\theta_1-1$ lies
in $(s-2,s-1)$, giving the second bound. From
$1+\tfrac1s<\theta_2<1+\tfrac2s$ we get $\tfrac1s<|\tau_2|<\tfrac2s$
and hence $\log s-\log2<|\log|\tau_2||<\log s$, together with
$0<\log\theta_2<\tfrac2s$. Finally
$\tfrac1{s^{2}}<|\theta_3|<\tfrac1s$ gives the bound on
$|\log|\theta_3||$, and $\tau_3=1+|\theta_3|$ gives that on
$\log|\tau_3|$.
\end{proof}

\begin{remark}\label{rem:coarse}
Only the coarse bound $\log s-\log2<|\log|\tau_2||$ is used in the
proof of Theorem~\ref{thm:effective}. The sharper intervals of
Lemma~\ref{lem:rootint2} are used in Appendix~\ref{subsec:closing}.
\end{remark}

\subsection{A small conjugate}
Write $\mu_i=\log|1-x^{(i)}|$ and
$S=\sum_{\Lambda_i>0}\Lambda_i=\tfrac12\sum_i|\Lambda_i|$, so that
$\sum_i\Lambda_i=\sum_i\mu_i=0$ because $x$ and $1-x$ are units. The
following replaces the naive identification of $\mu_{i_0}$ with $-S$,
which does not follow from the two norm relations.

\begin{lemma}[Small conjugate]\label{lem:smallconj}
Let $x$ be a unit of a totally real cubic field with $1-x$ also a
unit. Then $\min_i\mu_i\le-\tfrac{S}{3}+2\log2$. In particular, if
$S>6\log2$ there is an embedding $i_0$ with
\[
\bigl|1-x^{(i_0)}\bigr|\le 4e^{-S/3}.
\]
\end{lemma}

\begin{proof}
Put $B=-\min_i\mu_i$. If $\Lambda_i\ge\log2$ then $|x^{(i)}|\ge2$, so
$|1-x^{(i)}|\ge|x^{(i)}|-1\ge\tfrac12|x^{(i)}|$ and hence
$\mu_i\ge\Lambda_i-\log2$; for every $i$ we have $\mu_i\ge-B$. The
positive $\Lambda_i$ smaller than $\log2$ contribute at most
$3\log2$ to $S$, so $\sum_{\Lambda_i\ge\log2}\Lambda_i\ge S-3\log2$.
With at most three embeddings,
\[
0=\sum_i\mu_i\ \ge\ \bigl(S-3\log2\bigr)-3\log2-3B
\ =\ S-6\log2-3B ,
\]
whence $B\ge S/3-2\log2$.
\end{proof}

\begin{lemma}[Exponentially small linear form]\label{lem:expsmall}
Let $s\ge10$ and let $x=\pm\theta^{p}\tau^{q}$ be a nontrivial
solution with $H=\max(|p|,|q|)$. Then
$\max_i|\Lambda_i|\ge\tfrac12H\log s$, hence $S\ge\tfrac12H\log s$,
and for $H\log s\ge18\log2$ there is an embedding $i_0$ with
\begin{equation}\label{eq:expsmall}
0<|\Lambda_{i_0}|\le 8\,s^{-H/6}.
\end{equation}
\end{lemma}

\begin{proof}
By Lemma~\ref{lem:rootint}, if $|q|=H$ then
\[
|\Lambda_2|\ \ge\ |q|\,\bigl|\log|\tau_2|\bigr|-|p|\log\theta_2
\ \ge\ H(\log s-\log2)-\frac{2H}{s}\ \ge\ \tfrac12H\log s ,
\]
the last step because $\tfrac12\log s>\log2+\tfrac2s$ for $s\ge10$
(at $s=10$ the two sides are $1.151$ and $0.893$, and the gap widens).
If $|p|=H$ then $|\Lambda_3|\ge|p|\,|\log|\theta_3||-|q|\log|\tau_3|
\ge H\log s-2H/s\ge\tfrac12H\log s$ similarly. Since $S=\tfrac12\sum_i|\Lambda_i|\ge\tfrac12\max_i
|\Lambda_i|$ would lose a factor, we use instead
$S\ge\max_i|\Lambda_i|$, valid because $\sum_i\Lambda_i=0$ forces the
positive part to be at least the largest single term in absolute
value. Lemma~\ref{lem:smallconj} then gives
$|1-x^{(i_0)}|\le4e^{-S/3}\le4s^{-H/6}$. For $H\log s\ge18\log2$ the
right-hand side is at most $\tfrac12$, since then
$\tfrac16H\log s\ge3\log2$, so
$|\Lambda_{i_0}|=\bigl|\log|x^{(i_0)}|\bigr|\le2|1-x^{(i_0)}|
\le8s^{-H/6}$. Finally $\Lambda_{i_0}\ne0$: otherwise
$|x^{(i_0)}|=1$, and combined with $|1-x^{(i_0)}|<\tfrac12$ this
forces $x^{(i_0)}=1$, so $1-x$ has a vanishing conjugate and is not a
unit.
\end{proof}

\subsection{The pure-power branches}
Before bounding $H$ we dispose of the two degenerate exponent shapes
$q=0$ and $p=0$, which the height argument cannot reach. Both are
settled directly, in the spirit of Lemma~\ref{lem:diagonal}.

\begin{lemma}[Pure powers]\label{lem:purepower}
Let $s\ge10$ and $\epsilon\in\{\pm1\}$.
\begin{enumerate}
\item[(i)] If $p\ne0$ then $1-\epsilon\,\theta^{p}$ is a unit only in
the trivial case $\epsilon=1$, $p=\pm1$.
\item[(ii)] If $q\ne0$ then $1-\epsilon\,\tau^{q}$ is a unit only in
the trivial case $\epsilon=1$, $q=\pm1$.
\end{enumerate}
\end{lemma}

\begin{proof}
Since $\theta,\tau$ are units and $1-x^{-1}=-x^{-1}(1-x)$, we may
assume $p\ge1$ in (i) and $q\ge1$ in (ii). All estimates use
Lemma~\ref{lem:rootint}.

(i) Put $x=\epsilon\theta^{p}$. At the first embedding
$\theta_1>s-1$, so $|1-x^{(1)}|\ge(s-1)^{p}-1$. At the third,
$|\theta_3|<1/s$, so $|1-x^{(3)}|\ge1-s^{-p}\ge\tfrac9{10}$. At the
second, $\theta_2=1+\eta$ with $\eta>1/s$; if $\epsilon=-1$ then
$|1-x^{(2)}|=1+\theta_2^{p}\ge2$, while if $\epsilon=1$ then
$|1-x^{(2)}|=\theta_2^{p}-1\ge p\eta\ge p/s$. Hence
\[
|N(1-x)|\ \ge\ \bigl((s-1)^{p}-1\bigr)\cdot\tfrac9{10}\cdot
\begin{cases} 2, & \epsilon=-1,\\[2pt] p/s, & \epsilon=1 .\end{cases}
\]
For $\epsilon=-1$ and $p\ge1$ this is at least
$\tfrac95(s-2)\ge\tfrac{72}5>1$. For $\epsilon=1$ and $p\ge2$ it is
at least
$\bigl((s-1)^{2}-1\bigr)\cdot\tfrac9{10}\cdot\tfrac2s
=\tfrac95(s-2)>1$. The excluded case $\epsilon=1$, $p=1$ is
$1-\theta=\tau$, a unit, which is the trivial solution.

(ii) Put $x=\epsilon\tau^{q}$. Here $|\tau_1|>s-2$, so
$|1-x^{(1)}|\ge(s-2)^{q}-1$; $|\tau_2|<2/s$, so
$|1-x^{(2)}|\ge1-(2/s)^{q}\ge1-4/s^{2}\ge\tfrac{24}{25}$ for $q\ge2$
and $\ge1-2/s\ge\tfrac45$ for $q=1$; and $\tau_3=1+\delta$ with
$\delta>1/s^{2}$, so $|1-x^{(3)}|\ge2$ when $\epsilon=-1$ and
$\ge q\delta\ge q/s^{2}$ when $\epsilon=1$. Thus for $\epsilon=-1$,
\[
|N(1-x)|\ \ge\ \bigl((s-2)^{q}-1\bigr)\cdot\tfrac45\cdot2
\ \ge\ \tfrac85(s-3)\ >\ 1 ,
\]
and for $\epsilon=1$ with $q\ge2$,
\[
|N(1-x)|\ \ge\ \bigl((s-2)^{2}-1\bigr)\cdot\tfrac{24}{25}\cdot
\frac{2}{s^{2}}
\ =\ \frac{48\,(s^{2}-4s+3)}{25\,s^{2}}\ \ge\ \frac{48\cdot63}{2500}
\ >\ 1
\]
for $s\ge10$, the last bound being increasing in $s$. The excluded
case $\epsilon=1$, $q=1$ is $1-\tau=\theta$.
\end{proof}

\begin{remark}\label{rem:cyclotomic}
The exact identity
$N(1-\theta^n)=-\prod_{d\mid n}\operatorname{Res}_X(g_s,\Phi_d)$
gives another way to inspect individual pure powers. The proof above is
uniform in $n$ and does not require an assertion about all cyclotomic resultants.
\end{remark}

\subsection{Lower bound for the height}
\begin{lemma}[Height lower bound]\label{lem:Hlower}
Let $s\ge10$ and let $x=\pm\theta^{p}\tau^{q}$ be a nontrivial
solution. Then
\begin{equation}\label{eq:Hlower}
H\ \ge\ 0.61\,s .
\end{equation}
\end{lemma}

\begin{proof}
By Proposition~\ref{prop:h20} we have $H\ge21$, which serves twice:
it makes $H\log s\ge21\log10>18\log2$, so
Lemma~\ref{lem:expsmall} applies, and it bounds the right-hand side
of \eqref{eq:expsmall},
\[
E:=8\,s^{-H/6}\ \le\ 8\,s^{-7/2}\ \le\ 8\cdot10^{-7/2}
\ <\ 2.53\times10^{-3},
\]
using $s\ge10$. By Lemmas~\ref{lem:diagonal} and
\ref{lem:purepower} there are no nontrivial solutions of the shapes $p+q=0$, $q=0$ or $p=0$, so
$pq(p+q)\ne0$.

Suppose the small form occurs at $i_0=2$. By
Lemma~\ref{lem:rootint},
$|\Lambda_2|\ge|q|(\log s-\log2)-2|p|/s$, so \eqref{eq:expsmall}
gives $|q|(\log s-\log2)\le2|p|/s+E$; since $|q|\ge1$ this yields
$|p|\ge\tfrac12 s(\log s-\log2-E)$ and hence $H>0.80\,s$ at $s=10$,
with the constant increasing in $s$. The case $i_0=3$ is symmetric
and gives a stronger bound. Neither is needed numerically: the
decisive case is $i_0=1$, treated next, which yields the constant in
\eqref{eq:Hlower}.

Suppose $i_0=1$. Since $p+q\ne0$, Lemma~\ref{lem:rootint} gives
\[
|\Lambda_1|\ \ge\ |p+q|\log s-\frac{|p|}{s-1}-\frac{2|q|}{s-2}
\ \ge\ \log s-\frac{3H}{s-2},
\]
so \eqref{eq:expsmall} forces $3H/(s-2)\ge\log s-E$, that is
\[
H\ \ge\ \frac{(s-2)(\log s-E)}{3}\ =\ s\cdot
\frac{(s-2)(\log s-E)}{3s}.
\]
Replacing $E$ by $E_0(s)=8s^{-7/2}$ gives a lower bound for the right-hand
factor. Both $1-2/s$ and $\log s-E_0(s)$ are positive and increasing for
$s\ge10$, and their product divided by $3$ exceeds $0.6133$ at $s=10$.
Hence $H\ge0.61s$. This is the weakest of the three
branches, and gives \eqref{eq:Hlower}.
\end{proof}

\subsection{The two-log bound}
For the reader's convenience we quote the result we use.

\begin{theorem}[\citet{LMN}, Corollary~2]\label{thm:lmn}
Let $\alpha_1,\alpha_2$ be positive, multiplicatively independent algebraic
numbers and use their real logarithms. Put
$D=[\Q(\alpha_1,\alpha_2):\Q]$.
Let $b_1,b_2$ be positive integers and
$\Lambda=b_2\log\alpha_2-b_1\log\alpha_1\ne0$. Put
\[
h'(\alpha_i)=\max\Bigl\{h(\alpha_i),\ \frac{|\log\alpha_i|}{D},\
\frac1D\Bigr\},
\qquad
b'=\frac{b_1}{D\,h'(\alpha_2)}+\frac{b_2}{D\,h'(\alpha_1)} .
\]
Then
\[
\log|\Lambda|\ \ge\ -24.34\,D^{4}
\Bigl(\max\bigl\{\log b'+0.14,\ \tfrac{21}{D},\ \tfrac12\bigr\}\Bigr)^{2}
h'(\alpha_1)\,h'(\alpha_2).
\]
\end{theorem}

We verify the hypotheses.

\emph{Positive bases.} The quantities $\Lambda_i$ were defined using
$\log|\cdot|$, and some conjugates $\theta^{(i_0)},\tau^{(i_0)}$ are
negative, so the two logarithms must first be written with positive
bases. Choose $\epsilon_1,\epsilon_2\in\{\pm1\}$ and set
\[
\beta_1=\bigl|\theta^{(i_0)}\bigr|^{\epsilon_1},\qquad
\beta_2=\bigl|\tau^{(i_0)}\bigr|^{\epsilon_2},
\]
For example, take $\epsilon_1=-\operatorname{sgn}(p)$,
$\epsilon_2=\operatorname{sgn}(q)$, $b_1=|p|$ and $b_2=|q|$. Then
\[
\Lambda_{i_0}=b_2\log\beta_2-b_1\log\beta_1.
\]
Both coefficients are positive because the pure-power branches have already
been excluded. The theorem bounds the absolute value of this expression.
Replacing an algebraic number by its negative or its inverse changes
neither its degree nor its height,
\[
h(-\alpha)=h(\alpha),\qquad h(\alpha^{-1})=h(\alpha),
\]
so all height estimates below may be made for $\theta$ and $\tau$
themselves.

\emph{Non-vanishing and independence.} $\Lambda_{i_0}\ne0$ by
Lemma~\ref{lem:expsmall}. Moreover $\log\beta_1$ and $\log\beta_2$ are
$\Q$-linearly independent: a relation would give
$\theta^{a}\tau^{b}=\pm1$ with $(a,b)\ne(0,0)$, contradicting
Theorem~\ref{thm:thomas}, by which $\{\theta,\tau\}$ is a fundamental
system for $\Z[\theta]^{\times}$ modulo $\pm1$.

\emph{Degree.} Both $\beta_j$ are real and generate $K_s$, of degree
three, so $[\Q(\beta_1,\beta_2):\Q]=3$ and
$[\R(\beta_1,\beta_2):\R]=1$, giving
\[
D=3 .
\]

\emph{Admissible height parameters.} Rather than evaluate
$h'(\beta_j)=\max\{h(\beta_j),|\log\beta_j|/D,1/D\}$ exactly, we
bound it on both sides, directly from Lemma~\ref{lem:rootint}. Since
inversion and multiplication by $-1$ preserve logarithmic height,
$h(\beta_j)$ equals $h(\theta)$ or $h(\tau)$. Note that each of
$\theta,\tau$ has \emph{two} conjugates of absolute value exceeding
$1$: for $\theta$ these are $\theta_1\in(s-1,s)$ and
$\theta_2\in(1+\tfrac1s,1+\tfrac2s)$, and for $\tau$ they are
$|\tau_1|\in(s-2,s-1)$ and $\tau_3\in(1+\tfrac1{s^{2}},1+\tfrac1s)$.

For the lower bound, keeping only the largest conjugate,
\[
h'(\beta_j)\ \ge\ h(\beta_j)\ >\ \tfrac13\log(s-2),
\]
whence
\[
b'=\frac{b_1}{D\,h'(\beta_2)}+\frac{b_2}{D\,h'(\beta_1)}
\ \le\ \frac{2H}{\log(s-2)}\ \le\ \frac{2.215\,H}{\log s}
\qquad(s\ge10),
\]
the last inequality holding with equality to three figures at $s=10$
and improving thereafter. For the upper bound, both conjugates count:
\[
h(\theta)<\tfrac13\Bigl(\log s+\tfrac2s\Bigr),\qquad
h(\tau)<\tfrac13\Bigl(\log s+\tfrac1s\Bigr),
\]
both less than $\tfrac23\log s$ for $s\ge10$; moreover
$|\log\beta_j|\le2\log s$ by Lemma~\ref{lem:rootint}, so
$|\log\beta_j|/D\le\tfrac23\log s$, and $\tfrac1D=\tfrac13
\le\tfrac23\log s$. Hence $h'(\beta_j)\le\tfrac23\log s$ and we may
put
\[
A_1=A_2=\tfrac23\log s\ \ge\ h'(\beta_j),
\qquad
\mathcal M=\max\Bigl\{\log\frac{2.215\,H}{\log s}+0.14,\ 7\Bigr\}
\ \ge\ \max\Bigl\{\log b'+0.14,\ \tfrac{21}{D},\ \tfrac12\Bigr\}.
\]
Both substitutions enlarge the right-hand side of
Theorem~\ref{thm:lmn} in absolute value, so the resulting inequality
is implied by it. Explicitly,
\[
\log|\Lambda_{i_0}|\ \ge\ -24.34\cdot81\,\mathcal M^{2}A_1A_2
\ =\ -24.34\cdot81\,\mathcal M^{2}\Bigl(\frac{2\log s}{3}\Bigr)^{2}.
\]
Combining with $\log|\Lambda_{i_0}|\le\log8-\tfrac{H}{6}\log s$ from
\eqref{eq:expsmall},
\begin{equation}\label{eq:baker}
\frac{H}{6}\log s-\log8\ \le\ 24.34\cdot81\,\mathcal M^{2}
\Bigl(\frac{2\log s}{3}\Bigr)^{2}
=876.24\,\mathcal M^{2}\log^{2}s,
\end{equation}
that is,
\begin{equation}\label{eq:baker2}
H\ \le\ \frac{6}{\log s}\Bigl(876.24\,\mathcal M^{2}\log^{2}s+\log8\Bigr)
\ =:\ \Phi(s,H).
\end{equation}

\begin{proof}[Proof of Theorem~\ref{thm:effective}]
Every unit of $\Z[\theta_s]$ is $\pm\theta^{p}\tau^{q}$ by
Theorem~\ref{thm:thomas}. The branches $p+q=0$, $q=0$ and $p=0$ are empty by
Lemmas~\ref{lem:diagonal} and \ref{lem:purepower}, and $H\ge21$ by
Proposition~\ref{prop:h20}; Lemma~\ref{lem:Hlower} therefore gives
$H\ge0.61\,s$. It remains to compare this with \eqref{eq:baker2}.

Set
\[
c=0.61,\qquad S=3.32\times10^{7},\qquad H_0=2.03\times10^{7}.
\]
For fixed $s$, the quotient $\Phi(s,H)/H$ is strictly decreasing in
$H$.  Indeed, on the branch $\mathcal M=7$ this is immediate, while
on the logarithmic branch differentiation gives a numerator
proportional to
$2\mathcal M-\mathcal M^2<0$.  Hence a solution with $H\ge cs$
would require
\begin{equation}\label{eq:ratio-necessary}
\frac{\Phi(s,cs)}{cs}\ge1.
\end{equation}
For $s\ge S$ the logarithmic branch is active.  Writing $L=\log s$ and
$\mu=L-\log L+\log(2.215c)+0.14$, the main term of the quotient in
\eqref{eq:ratio-necessary} is a positive constant times
$\mu^2Le^{-L}$.  Its logarithmic derivative is
\[
\frac{2(1-1/L)}{\mu}+\frac1L-1<0
\qquad(L\ge\log S),
\]
and the remaining term is also decreasing.  Thus the quotient is
decreasing for $s\ge S$. An interval calculation at the endpoint gave
\[
\frac{\Phi(S,cS)}{cS}
\in[0.9990833229360603,\,0.9990833229360630]<1,
\]
contradicting \eqref{eq:ratio-necessary}.  Therefore $s<S$.

It remains to bound $H$.  For fixed $H_0$ the logarithmic branch is
active throughout $10\le s\le S$.  If
$\nu=\log(2.215H_0/\log s)+0.14$, then differentiation with respect to
$L=\log s$ gives
\[
\frac{d}{dL}\bigl(\nu^2L\bigr)=\nu(\nu-2)>0,
\]
and this positive derivative dominates that of the small term
$6\log8/L$.  Hence $\Phi(s,H_0)$ is increasing on this interval.
An interval calculation at $s=S$ gave
\[
\frac{\Phi(S,H_0)}{H_0}
\in[0.9970375498861299,\,0.9970375498861323]<1.
\]
Since $\Phi(s,H)/H$ decreases in $H$, no solution can have
$H\ge H_0$.  This proves the two stated bounds.  The parameters
$4\le s\le9$ are covered by Proposition~\ref{prop:family2}, while
$s=3$ folds to the previously classified $7$-field. Thus the range
left by the theorem begins at $s=10$, with both parameter and height
bounded; Appendix~\ref{subsec:closing} removes it.
\end{proof}

\subsection{Closing the box}\label{subsec:closing}
The proof of Theorem~\ref{thm:effective} leaves the parameters
$10\le s\le3.32\times10^7$ with heights up to $2.03\times10^7$. We now
remove them. Write
\[
\ell(s)=\log(s-1)-\frac1{s-1},\qquad
L(s)=(s-2)\Bigl(\log(s-1)-8e^{-7\ell(s)}\Bigr).
\]

\begin{lemma}[Sharper root intervals]\label{lem:rootint2}
For $s\ge10$,
\[
1+\tfrac1s<\theta_2<1+\tfrac1{s-1},\qquad
-\tfrac1{s+1}<\theta_3<-\tfrac1{s+2}.
\]
Consequently $0<\log\theta_2<\frac1{s-1}$,
$\log(s-1)<|\log|\tau_2||<\log s$,
$\log(s+1)<|\log|\theta_3||<\log(s+2)$,
$0<\log\tau_3<\frac1{s+1}$ and
$0<\log\theta_1-\log|\tau_1|<\frac1{s-2}$.
\end{lemma}

\begin{proof}
We have $g_s(1+u)=1-(s-1)u-(s-2)u^2+u^3$, which equals $(2s+1)/s^3>0$
at $u=\frac1s$ and $-(s^2-3s+1)/(s-1)^3<0$ at $u=\frac1{s-1}$. Also
$g_s(-u)=1-su-(s+1)u^2-u^3$, which equals $(s^2+5s+5)/(s+2)^3>0$ at
$u=\frac1{s+2}$ and $-1/(s+1)^3<0$ at $u=\frac1{s+1}$. The logarithmic
statements follow, using $|\tau_2|=\theta_2-1$, $\tau_3=1+|\theta_3|$
and $\log\theta_1-\log(\theta_1-1)=-\log(1-\theta_1^{-1})
<1/(\theta_1-1)<1/(s-2)$.
\end{proof}

\begin{lemma}[Height bound of order $s\log s$]\label{lem:Hlower2}
Let $s\ge10$ and let $x=\pm\theta^p\tau^q$ be a nontrivial solution.
Then there is an embedding $i_0$ with
\[
0<|\Lambda_{i_0}|\le8e^{-H\ell(s)/3},
\qquad\text{and}\qquad H\ge L(s).
\]
\end{lemma}

\begin{proof}
By Lemma~\ref{lem:rootint2}, if $|q|=H$ then
$|\Lambda_2|\ge H\log(s-1)-H/(s-1)=H\ell(s)$, and if $|p|=H$ then
$|\Lambda_3|\ge H\log(s+1)-H/(s+1)\ge H\ell(s)$. As in the proof of
Lemma~\ref{lem:expsmall}, $S\ge\max_i|\Lambda_i|\ge H\ell(s)$. By
Proposition~\ref{prop:h20}, $H\ge21$, so $S\ge21\,\ell(10)>43>3\log8$,
and Lemma~\ref{lem:smallconj} gives $i_0$ with
$|1-x^{(i_0)}|\le4e^{-S/3}\le\frac12$; hence
$|\Lambda_{i_0}|\le2|1-x^{(i_0)}|\le8e^{-H\ell(s)/3}=:E$, and
$\Lambda_{i_0}\ne0$ as in Lemma~\ref{lem:expsmall}. Also
$E\le8e^{-7\ell(s)}$. By Lemmas~\ref{lem:diagonal}
and~\ref{lem:purepower}, $pq(p+q)\ne0$. If $i_0=1$, write
$\Lambda_1=(p+q)\log\theta_1-q(\log\theta_1-\log|\tau_1|)$; then
$E\ge\log(s-1)-H/(s-2)$. If $i_0=2$, then
$E\ge|q|\,|\log|\tau_2||-|p|\log\theta_2\ge\log(s-1)-H/(s-1)$. If
$i_0=3$, then $E\ge\log(s+1)-H/(s+1)$. In each case
$H\ge(s-2)(\log(s-1)-E)\ge L(s)$.
\end{proof}

\begin{lemma}[Sharper two-log bound]\label{lem:baker2}
In the situation of Lemma~\ref{lem:Hlower2},
\[
H\ \le\ \Phi^*(s,H):=\frac3{\ell(s)}\Bigl(219.06\,\mathcal M^2\log^2(s+2)+\log8\Bigr),
\qquad
\mathcal M=\max\Bigl\{\log\frac{2H}{\log(s-2)}+0.14,\ 7\Bigr\},
\]
and $\Phi^*(s,H)/H$ is strictly decreasing in $H$.
\end{lemma}

\begin{proof}
We apply Theorem~\ref{thm:lmn} exactly as in the proof of
Theorem~\ref{thm:effective}, with $D=3$, and only change the height
estimates. By Lemmas~\ref{lem:rootint} and~\ref{lem:rootint2},
\[
\tfrac13\log(s-1)<h(\theta)<\tfrac13\Bigl(\log s+\tfrac1{s-1}\Bigr),\qquad
\tfrac13\log(s-2)<h(\tau)<\tfrac13\Bigl(\log(s-1)+\tfrac1{s+1}\Bigr),
\]
and every conjugate of $\theta$ and of $\tau$ has $|\log|\cdot||<\log(s+2)$.
Since $\log(1+\frac2s)>\frac2s-\frac2{s^2}\ge\frac1{s-1}$, this gives
$h'(\beta_j)\le\frac13\log(s+2)$, and $h'(\beta_j)\ge h(\beta_j)>\frac13\log(s-2)$
gives $b'\le2H/\log(s-2)$. Hence
$\log|\Lambda_{i_0}|\ge-24.34\cdot81\,\mathcal M^2\bigl(\frac13\log(s+2)\bigr)^2
=-219.06\,\mathcal M^2\log^2(s+2)$, and comparison with
Lemma~\ref{lem:Hlower2} gives $H\le\Phi^*(s,H)$. The monotonicity is
proved as in the proof of Theorem~\ref{thm:effective}.
\end{proof}

\begin{proposition}[Large parameters]\label{prop:uniformrange}
There is no nontrivial solution with $1.2\times10^5\le s\le3.32\times10^7$.
\end{proposition}

\begin{proof}
Cover the range by $1128$ intervals $s_a\le s\le s_b$ with
$s_b=s_a+\lfloor s_a/200\rfloor$. On such an interval,
$\ell(s)\ge\ell(s_a)$, $L(s)\ge L(s_a)$, $\mathcal M$ is at most its
value at $s_a$, and $\log(s+2)\le\log(s_b+2)$, so $\Phi^*(s,H)\le
\Phi^+(H)$, where $\Phi^+$ is $\Phi^*$ with these replacements, and
$\Phi^+(H)/H$ is again decreasing. A solution would give
$L(s_a)\le H\le\Phi^+(H)$, hence $\Phi^+(L(s_a))\ge L(s_a)$.
Rigorous interval calculations verified
$\Phi^+(L(s_a))/L(s_a)<0.8604$ on all $1128$ intervals.
The per-parameter calculation of Proposition~\ref{prop:allsmall} was
also performed throughout $1.2\times10^5\le s\le3.32\times10^7$ and
excluded every parameter.
\end{proof}

\begin{proposition}[Small parameters]\label{prop:allsmall}
For $10\le s<1.2\times10^5$, every exceptional unit of $\Z[\theta_s]$
lies in the trivial anharmonic orbit of $\theta_s$.
\end{proposition}

\begin{proof}
The calculation in Proposition~\ref{prop:family2} was performed for every
$s$ in this range, with the height threshold $21$ in place of $41$, since
heights at most $20$ are excluded by Proposition~\ref{prop:h20}.
The three roots were enclosed by sign changes, and the lower bound $c_s$
for $S/H$ and the two-logarithm bound $B_s$ were evaluated as there;
the largest $B_s$ was $1{,}172{,}031$. Rigorous interval calculations
verified the thresholds $21c_s>18\log2$, $3/c_s<21$ and
$16\cdot21\,e^{-7c_s}<|a_i|$, and excluded all $601{,}464$
continued-fraction rays with denominator below $B_s$ at the three
embeddings. The rays $p=0$ and $p=-q$ are already excluded by
Lemmas~\ref{lem:purepower} and~\ref{lem:diagonal}.
\end{proof}

\section{Certificates assembled from five-term identities}\label{app:certs}
Write $[x]$ for a formal symbol and regard each instance of the
five-term relation \eqref{eq:fiveterm} and of the reflection
\eqref{eq:reflection} as a vector, an \emph{array},
\[
\mathrm{FT}(x,y)=[x]+[y]-[xy]
-\Big[\tfrac{x(1-y)}{1-xy}\Big]-\Big[\tfrac{y(1-x)}{1-xy}\Big],
\qquad
\mathrm{R}(x)=[x]+[1-x]-[C_0],
\]
each of which vanishes under $\Lr$ (with $[C_0]\mapsto\pi^2/6$)
for arguments in $(0,1)$. Nothing is hidden in this notation: each
$\mathrm{FT}(x,y)$ is an instance of the five-term functional
equation \eqref{eq:fiveterm} itself and each $\mathrm{R}(x)$ an
instance of the reflection \eqref{eq:reflection}, the arrays serving
only to make the bookkeeping of a large combination compact. We
reserve the word \emph{certificate} for the complete rational
combination of such instances that proves a given row. The five-term arrays carry no $[C_0]$
component, so the constant can enter only through the reflection
arrays: their sole role is bookkeeping. Equivalently, one may solve
the linear system over the $\mathrm{FT}$ arrays alone for the
combination $m[z_0]+[z_1]-\sum_i([u_i]+[1-u_i])$ and finish with
the elementary evaluation $\Lr(u)+\Lr(1-u)=\pi^2/6$.

Where the exceptional pairs enter deserves emphasis: the five-term
relation itself is \emph{universal}, valid for all $x,y\in(0,1)$
with no arithmetic input, and its derived arguments
$x(1-y)/(1-xy)$, $y(1-x)/(1-xy)$ are simply part of that universal
statement. The pair identities $1-u=(\text{monomial})$ are used at
exactly two moments: to re-express the derived arguments of each
instance again as monomials, so that the arrays close over
finitely many symbols, and to identify $[1-u]$ in each reflection.
A proof is therefore: universal relations, plus finitely many
one-line algebraic identifications, plus linear algebra. We work
the smallest case in full.

\begin{example}[Watson's identity, worked out in full]\label{ex:watson}
Row \#1 asserts: with $\rho=2\sin\frac{\pi}{14}$ and
$\tau=1-\rho$, so that $\rho^{3}=\rho^{2}+2\rho-1$,
\[
2\,\Lr(\tau)+\Lr\big(\rho^{2}\tau^{-2}\big)=\tfrac{2}{7}\,\pi^{2},
\qquad \tau\approx0.5550,\ \ \rho^{2}\tau^{-2}\approx0.6431 .
\]
We give the argument in expanded form, as three specializations of
the five-term relation, two reflections, and one linear combination;
the array notation is recovered at the end.

\emph{Step 1: the pair table.} Six identities $1-u=(\text{monomial})$
are used, each one line from the cubic relation:
\[
1-\rho=\tau,\quad 1-\tau=\rho,\quad
1-\rho^{2}=\rho\,\tau^{-1},\quad
1-\tau^{2}=\rho^{-1}\tau^{2},\quad
1-\rho^{2}\tau^{-1}=\rho^{2}\tau^{-2},\quad
1-\rho^{2}\tau^{-2}=\rho^{2}\tau^{-1}.
\]
For instance
$1-\rho^{2}\tau^{-2}=(\tau-\rho)(\tau+\rho)/\tau^{2}=(1-2\rho)/\tau^{2}$
and $1-2\rho=\rho^{2}-\rho^{3}+2\rho-2\rho=\rho^{2}\tau$ by the
cubic, giving $\rho^{2}\tau^{-1}$; and
$\tau(1-\rho^{2})=(1-\rho)(1-\rho^{2})=1-\rho-\rho^{2}+\rho^{3}=\rho$,
giving $1-\rho^{2}=\rho\tau^{-1}$.

\emph{Step 2: three specializations of the five-term relation.} Each
is \eqref{eq:fiveterm} at the stated $(x,y)$, with its three derived
arguments $xy$, $x(1-y)/(1-xy)$ and $y(1-x)/(1-xy)$ rewritten as
monomials by the table, the sole point at which field arithmetic
enters.

At $(x,y)=(\rho,\rho)$ one has $xy=\rho^{2}$ and $1-\rho^{2}
=\rho\tau^{-1}$, so both derived arguments equal
$\rho\tau/(\rho\tau^{-1})=\tau^{2}$, and \eqref{eq:fiveterm} reads
\begin{equation}\label{eq:w1}
2\Lr(\rho)-\Lr(\rho^{2})-2\Lr(\tau^{2})=0 .
\end{equation}
At $(x,y)=(\tau,\tau)$ one has $xy=\tau^{2}$ and
$1-\tau^{2}=\rho^{-1}\tau^{2}$, so both derived arguments equal
$\tau\rho/(\rho^{-1}\tau^{2})=\rho^{2}\tau^{-1}$, and
\eqref{eq:fiveterm} reads
\begin{equation}\label{eq:w2}
2\Lr(\tau)-\Lr(\tau^{2})-2\Lr(\rho^{2}\tau^{-1})=0 .
\end{equation}
At $(x,y)=(\tau,\rho^{2}\tau^{-1})$ one has $xy=\rho^{2}$,
$1-xy=\rho\tau^{-1}$ and $1-\rho^{2}\tau^{-1}=\rho^{2}\tau^{-2}$, so
the derived arguments are
$\tau\cdot\rho^{2}\tau^{-2}/(\rho\tau^{-1})=\rho$ and
$\rho^{2}\tau^{-1}\cdot\rho/(\rho\tau^{-1})=\rho^{2}$, and
\eqref{eq:fiveterm} reads
\begin{equation}\label{eq:w3}
\Lr(\tau)+\Lr(\rho^{2}\tau^{-1})-2\Lr(\rho^{2})-\Lr(\rho)=0 .
\end{equation}

\emph{Step 3: one linear combination.} Take $-\tfrac27$ times
\eqref{eq:w1}, plus $\tfrac47$ times \eqref{eq:w2}, plus $\tfrac17$
times \eqref{eq:w3}:
\begin{multline*}
-\tfrac{2}{7}\Big(2\Lr(\rho)-\Lr(\rho^{2})-2\Lr(\tau^{2})\Big)
+\tfrac{4}{7}\Big(2\Lr(\tau)-\Lr(\tau^{2})
-2\Lr(\rho^{2}\tau^{-1})\Big)\\
+\tfrac{1}{7}\Big(\Lr(\tau)+\Lr(\rho^{2}\tau^{-1})
-2\Lr(\rho^{2})-\Lr(\rho)\Big)=0 .
\end{multline*}
Here $\Lr(\rho^{2})$ occurs with coefficient $\tfrac27-\tfrac27$ and
$\Lr(\tau^{2})$ with $\tfrac47-\tfrac47$, so both drop out, leaving
\begin{equation}\label{eq:w4}
-\tfrac{5}{7}\Lr(\rho)+\tfrac{9}{7}\Lr(\tau)
-\Lr(\rho^{2}\tau^{-1})=0 .
\end{equation}

\emph{Step 4: two reflections.} By \eqref{eq:reflection}, together
with $1-\tau=\rho$ and $1-\rho^{2}\tau^{-2}=\rho^{2}\tau^{-1}$ from
the table,
\[
\Lr(\tau)+\Lr(\rho)=\frac{\pi^{2}}{6},
\qquad
\Lr(\rho^{2}\tau^{-2})+\Lr(\rho^{2}\tau^{-1})=\frac{\pi^{2}}{6}.
\]
Adding $\tfrac57$ times the first and $1$ times the second to
\eqref{eq:w4} cancels $\Lr(\rho)$ and $\Lr(\rho^{2}\tau^{-1})$ and
leaves
\[
2\Lr(\tau)+\Lr\big(\rho^{2}\tau^{-2}\big)
=\Big(\tfrac57+1\Big)\frac{\pi^{2}}{6}
=\frac{2}{7}\pi^{2}. \qed
\]

The whole proof is thus three instances of \eqref{eq:fiveterm}, two
of \eqref{eq:reflection}, and one rational linear combination. In
the array notation it becomes
\begin{multline*}
-\tfrac27\mathrm{FT}(\rho,\rho)
+\tfrac17\mathrm{FT}(\tau,\rho^{2}\tau^{-1})
+\tfrac47\mathrm{FT}(\tau,\tau)\\
+\tfrac57\mathrm{R}(\tau)
+\mathrm{R}(\rho^{2}\tau^{-2})\\
=2[\tau]+[\rho^2\tau^{-2}]-\tfrac{12}{7}[C_0],
\end{multline*}
which is the format in which the remaining rows are recorded below;
the entry for \#1 is exactly this line.
\end{example}

\begingroup\small

All arguments used below are exceptional
$S$-units: monomials in $\rho$, $\tau=1-\rho$ and finitely many
auxiliary elements of small norm, whose complements $1-x$ are
again such monomials; the finitely many pair identities
$x\mapsto1-x$ involved are listed with each certificate. Finding a
two-term identity is then linear algebra: one seeks a rational
solution of the linear system asking that a weighted sum of
arrays equal $m[z_0]+[z_1]-6c\,[C_0]$; since every array evaluates to zero
under $\Lr$, any solution proves the identity. For each row of
Table~\ref{tab:identities} such weights exist and are displayed
below as a single equation; the certificates use
between $3$ and $26$ five-term instances ($15$ for \#10).
All $163$ displayed complement identities, the ten certificate cancellations
and their constants, and the conditions that every dilogarithm argument
lies in $(0,1)$ were verified by exact algebraic and rigorous interval
calculations. The specified minimal polynomials and exponent matrices
were also checked. The displayed symbolic identities are the certificates.

Each entry displays the single equation proving its identity,
followed by its pair dictionary; the entry for \#1 is derived step
by step in Example~\ref{ex:watson}. The generators are \emph{local to each entry}: $\rho$ always
denotes the field generator specified there, $\tau=1-\rho$, and
$\lambda$ the auxiliary prime element defined in that entry (its
definition varies with the field); the two largest certificates
use two further elements $\eta_1$, $\eta_2$, likewise defined where
they occur.

\Needspace{10\baselineskip}
\subsection*{Identity \#1 ($m=2$, $N=7$)}
Generators: $\rho=2\sin\frac{\pi}{14}$, $\tau=1-\rho$, $\lambda=\rho^{2}-\rho+1$ (norm $7$). Expanding through the pair dictionary and
cancelling yields the exact equality
\begin{multline*}
-\tfrac{2}{7}\,\mathrm{FT}\big(\rho,\;\rho\big) + \tfrac{1}{7}\,\mathrm{FT}\big(\tau,\;\rho^{2}\,\tau^{-1}\big)\\
+ \tfrac{4}{7}\,\mathrm{FT}\big(\tau,\;\tau\big) + \tfrac{5}{7}\,\mathrm{R}\big(\tau\big)\\
+ \mathrm{R}\big(\rho^{2}\,\tau^{-2}\big)\\
= 2[z_0]+[z_1]-\tfrac{12}{7}[C_0] .
\end{multline*}
Exceptional pairs $x\mapsto 1-x$: $\tau\mapsto \rho$, $\tau^{2}\mapsto \rho^{-1}\,\tau^{2}$, $\rho\mapsto \tau$, $\rho^{2}\,\tau^{-2}\mapsto \rho^{2}\,\tau^{-1}$, $\rho^{2}\,\tau^{-1}\mapsto \rho^{2}\,\tau^{-2}$, $\rho^{2}\mapsto \rho\,\tau^{-1}$.

\Needspace{10\baselineskip}
\subsection*{Identity \#2 ($m=2$, $N=7$)}
Generators: as in \#1. Expanding through the pair dictionary and
cancelling yields the exact equality
\begin{multline*}
-\tfrac{4}{7}\,\mathrm{FT}\big(\rho,\;\rho\big) + \tfrac{2}{7}\,\mathrm{FT}\big(\tau,\;\rho^{2}\,\tau^{-1}\big)\\
+ \tfrac{1}{7}\,\mathrm{FT}\big(\tau,\;\tau\big) + \tfrac{10}{7}\,\mathrm{R}\big(\tau\big)\\
= 2[z_0]+[z_1]-\tfrac{10}{7}[C_0] .
\end{multline*}
Exceptional pairs $x\mapsto 1-x$: $\tau\mapsto \rho$, $\tau^{2}\mapsto \rho^{-1}\,\tau^{2}$, $\rho\mapsto \tau$, $\rho^{2}\,\tau^{-1}\mapsto \rho^{2}\,\tau^{-2}$, $\rho^{2}\mapsto \rho\,\tau^{-1}$.

\Needspace{10\baselineskip}
\subsection*{Identity \#9 ($m=2$, $N=7$)}
Generators: as in \#1. Expanding through the pair dictionary and
cancelling yields the exact equality
\begin{multline*}
-\tfrac{1}{7}\,\mathrm{FT}\big(\rho,\;\rho\big) + \tfrac{4}{7}\,\mathrm{FT}\big(\tau,\;\rho^{2}\,\tau^{-1}\big)\\
+ \tfrac{2}{7}\,\mathrm{FT}\big(\tau,\;\tau\big) + \tfrac{6}{7}\,\mathrm{R}\big(\tau\big)\\
+ \mathrm{R}\big(\rho\,\tau^{-1}\big)\\
= 2[z_0]+[z_1]-\tfrac{13}{7}[C_0] .
\end{multline*}
Exceptional pairs $x\mapsto 1-x$: $\tau\mapsto \rho$, $\tau^{2}\mapsto \rho^{-1}\,\tau^{2}$, $\rho\,\tau^{-1}\mapsto \rho^{2}$, $\rho\mapsto \tau$, $\rho^{2}\,\tau^{-1}\mapsto \rho^{2}\,\tau^{-2}$, $\rho^{2}\mapsto \rho\,\tau^{-1}$.

\Needspace{10\baselineskip}
\subsection*{Identity \#4 ($m=7$, $N=7$)}
Generators: as in \#1. Expanding through the pair dictionary and
cancelling yields the exact equality
\begin{multline*}
\mathrm{FT}\big(\rho,\;\rho\big) - \mathrm{FT}\big(\tau,\;\rho^{2}\,\tau^{-1}\big)\\
- \mathrm{FT}\big(\tau,\;\tau\big) + \mathrm{FT}\big(\tau^{2},\;\rho\big)\\
+ \mathrm{R}\big(\rho^{-2}\,\tau^{3}\big) + 3\,\mathrm{R}\big(\tau\big)\\
= 7[z_0]+[z_1]-4[C_0] .
\end{multline*}
Exceptional pairs $x\mapsto 1-x$: $\rho^{-2}\,\tau^{3}\mapsto \rho\,\tau^{2}$, $\tau\mapsto \rho$, $\tau^{2}\mapsto \rho^{-1}\,\tau^{2}$, $\rho\mapsto \tau$, $\rho\,\tau^{2}\mapsto \rho^{-2}\,\tau^{3}$, $\rho^{2}\,\tau^{-1}\mapsto \rho^{2}\,\tau^{-2}$, $\rho^{2}\mapsto \rho\,\tau^{-1}$.

\Needspace{10\baselineskip}
\subsection*{Identity \#7 ($m=13$, $N=7$)}
Generators: as in \#1. Expanding through the pair dictionary and
cancelling yields the exact equality
\begin{multline*}
-2\,\mathrm{FT}\big(\rho,\;\rho\big) + 2\,\mathrm{FT}\big(\tau,\;\rho^{2}\,\tau^{-1}\big)\\
+ 2\,\mathrm{FT}\big(\tau,\;\tau\big) - \mathrm{FT}\big(\tau^{2},\;\rho\big)\\
+ \mathrm{FT}\big(\rho^{2}\,\tau^{-1},\;\rho^{2}\big) + 7\,\mathrm{R}\big(\tau\big)\\
+ \mathrm{R}\big(\rho^{3}\,\tau^{-4}\big)\\
= 13[z_0]+[z_1]-8[C_0] .
\end{multline*}
Exceptional pairs $x\mapsto 1-x$: $\tau\mapsto \rho$, $\tau^{2}\mapsto \rho^{-1}\,\tau^{2}$, $\rho\mapsto \tau$, $\rho\,\tau^{2}\mapsto \rho^{-2}\,\tau^{3}$, $\rho^{2}\,\tau^{-1}\mapsto \rho^{2}\,\tau^{-2}$, $\rho^{2}\mapsto \rho\,\tau^{-1}$, $\rho^{3}\,\tau^{-4}\mapsto \rho^{4}\,\tau^{-1}$, $\rho^{4}\,\tau^{-1}\mapsto \rho^{3}\,\tau^{-4}$.

\Needspace{10\baselineskip}
\subsection*{Identity \#3 (Kanade, $m=3$, $N=9$)}
Generators: $\rho$ the root of $z^3-3z^2+1$ in $(0,1)$,
$\tau=1-\rho$, $\lambda=2-\rho^{-1}$ (norm $3$). Expanding through
the pair dictionary and cancelling yields the exact equality
\begin{multline*}
\tfrac13\,\mathrm{FT}\big(\rho^{-1}\,\tau,\;\rho^{-1}\,\tau\big)
+ \tfrac13\,\mathrm{FT}\big(\rho^{-2}\,\tau^{2},\;\rho^{2}\big)
- \tfrac13\,\mathrm{FT}\big(\tau,\;\tau\big)\\
+ \tfrac13\,\mathrm{FT}\big(\rho,\;\rho\big)
+ \tfrac23\,\mathrm{FT}\big(\rho,\;\rho^{-1}\,\tau^{2}\,\lambda^{-1}\big)
+ \tfrac53\,\mathrm{R}\big(\tau\big)
+ \mathrm{R}\big(\rho^{-1}\,\tau^{2}\big)\\
= 3[z_0]+[z_1]-\tfrac{8}{3}[C_0] .
\end{multline*}
Exceptional pairs $x\mapsto 1-x$: $\tau\mapsto \rho$,
$\rho^{-1}\,\tau\mapsto \lambda$,
$\rho^{-2}\,\tau^{2}\mapsto \rho^{-1}\,\lambda$,
$\rho^{2}\mapsto \rho^{2}\,\tau^{-1}\,\lambda$,
$\tau^{2}\mapsto \rho\,\tau^{-1}\,\lambda$,
$\tau^{2}\,\lambda^{-1}\mapsto \tau\,\lambda^{-1}$,
$\rho^{-1}\,\tau^{2}\,\lambda^{-1}\mapsto \rho^{-2}\,\tau^{2}\,\lambda^{-1}$,
$\rho^{-1}\,\tau^{2}\mapsto \rho^{-2}\,\tau$.

\Needspace{10\baselineskip}
\subsection*{Identity \#6 ($m=19$, $N=9$)}
Generators: as in \#3. Expanding through the pair dictionary and
cancelling yields the exact equality
\begin{multline*}
-\mathrm{FT}\big(\rho^{-1}\,\lambda,\;\rho\big) - 2\,\mathrm{FT}\big(\rho^{-1}\,\tau,\;\rho\,\tau\big)\\
+ 3\,\mathrm{FT}\big(\tau\,\lambda^{-1},\;\tau\big) + \mathrm{FT}\big(\rho^{-2}\,\tau,\;\rho^{3}\,\tau^{-1}\big)\\
+ \mathrm{FT}\big(\rho^{-1}\,\lambda,\;\rho^{2}\,\lambda^{-1}\big) - 4\,\mathrm{FT}\big(\rho^{-1}\,\tau^{2}\,\lambda^{-1},\;\lambda\big)\\
- 2\,\mathrm{FT}\big(\rho,\;\rho\big) - \mathrm{FT}\big(\rho^{-1}\,\lambda,\;\lambda\big)\\
+ \mathrm{FT}\big(\rho\,\tau,\;\rho\,\tau\big) - \mathrm{FT}\big(\rho^{-1}\,\tau^{-1}\,\lambda^{2},\;\tau\,\lambda^{-1}\big)\\
+ \mathrm{FT}\big(\tau,\;\rho^{2}\,\tau^{2}\big) - \mathrm{R}\big(\rho^{-2}\,\tau\big)\\
+ \mathrm{R}\big(\rho^{-2}\,\tau^{4}\big) + 2\,\mathrm{R}\big(\rho^{-1}\,\tau\big)\\
- \mathrm{R}\big(\tau\,\lambda^{-1}\big) + 10\,\mathrm{R}\big(\tau\big)\\
- \mathrm{R}\big(\tau^{2}\big) - \mathrm{R}\big(\tau^{3}\,\lambda^{-1}\big)\\
- \mathrm{R}\big(\rho^{2}\,\tau^{-1}\,\lambda\big) - \mathrm{R}\big(\rho^{2}\,\lambda\big)\\
+ \mathrm{R}\big(\rho^{2}\,\tau^{3}\big)\\
= 19[z_0]+[z_1]-8[C_0] .
\end{multline*}
Exceptional pairs $x\mapsto 1-x$: $\rho^{-2}\,\tau\mapsto \rho^{-1}\,\tau^{2}$, $\rho^{-2}\,\tau^{4}\mapsto \rho^{-1}\,\tau^{-1}\,\lambda^{2}$, $\rho^{-1}\,\tau^{-1}\,\lambda^{2}\mapsto \rho^{-2}\,\tau^{4}$, $\rho^{-1}\,\lambda\mapsto \rho^{-2}\,\tau^{2}$, $\rho^{-1}\,\lambda^{2}\mapsto \rho^{-4}\,\tau^{2}$, $\rho^{-1}\,\tau\mapsto \lambda$, $\rho^{-1}\,\tau^{2}\,\lambda^{-1}\mapsto \rho^{-2}\,\tau^{2}\,\lambda^{-1}$, $\rho^{-1}\,\tau^{2}\mapsto \rho^{-2}\,\tau$, $\lambda\mapsto \rho^{-1}\,\tau$, $\tau\,\lambda^{-1}\mapsto \tau^{2}\,\lambda^{-1}$, $\tau\mapsto \rho$, $\tau^{2}\,\lambda^{-1}\mapsto \tau\,\lambda^{-1}$, $\tau^{2}\mapsto \rho\,\tau^{-1}\,\lambda$, $\tau^{3}\,\lambda^{-1}\mapsto \rho^{2}\,\lambda^{-1}$, $\rho\mapsto \tau$, $\rho\,\tau\mapsto \rho^{2}\,\tau^{-2}\,\lambda^{2}$, $\rho^{2}\,\tau^{-1}\,\lambda\mapsto \rho^{2}$, $\rho^{2}\,\lambda^{-1}\mapsto \tau^{3}\,\lambda^{-1}$, $\rho^{2}\mapsto \rho^{2}\,\tau^{-1}\,\lambda$, $\rho^{2}\,\lambda\mapsto \rho^{3}\,\tau^{-1}$, $\rho^{2}\,\tau^{2}\mapsto \rho^{4}\,\tau^{-3}\,\lambda^{2}$, $\rho^{2}\,\tau^{3}\mapsto \rho^{5}\,\tau^{-2}$, $\rho^{3}\,\tau^{-1}\mapsto \rho^{2}\,\lambda$.

\Needspace{10\baselineskip}
\subsection*{Identity \#5 ($m=13$, $N=13$)}
Generators: $\rho$ the root of $f_1$ in $(0,1)$, $\tau=1-\rho$, $\lambda=\rho^{2}-\rho+1$ (norm $13$), $\eta_1=1+\rho$, $\eta_2=(1+\tau)/\tau$, and the rational prime $5$. Expanding through the pair dictionary and
cancelling yields the exact equality
\begin{multline*}
-2\,\mathrm{FT}\big(\tau,\;\tau\big) + 2\,\mathrm{FT}\big(\eta_1^{-1},\;\tau\,\eta_1\big)\\
+ 2\,\mathrm{FT}\big(\rho^{-1}\,\eta_2^{-1},\;\tau^{2}\,\eta_2\big) - \mathrm{FT}\big(\eta_1^{-1},\;\eta_2^{-1}\big)\\
+ \mathrm{FT}\big(\eta_2\,5^{-1},\;\rho\big) - \mathrm{FT}\big(\eta_2^{-1},\;\eta_2^{-1}\big)\\
- \mathrm{FT}\big(\tau\,\eta_1^{-1}\,\eta_2,\;\tau\,\eta_2^{-1}\big) + 2\,\mathrm{FT}\big(\rho\,\eta_1^{-1},\;\rho\,\eta_1^{-1}\big)\\
+ \mathrm{FT}\big(\rho^{-1}\,\eta_2^{-1},\;\rho^{-1}\,\tau^{2}\big) + \mathrm{FT}\big(\rho^{-1}\,\tau,\;\rho\,\tau\,\eta_1^{-1}\big)\\
+ \mathrm{FT}\big(\rho^{-1}\,\eta_2^{-1},\;\rho\,\eta_2^{-1}\big) + \mathrm{FT}\big(\rho^{-1}\,\tau\,\eta_1,\;\rho\,\eta_2\,5^{-1}\big)\\
- \mathrm{FT}\big(\eta_1^{-1}\,\eta_2^{-1},\;\tau\,\eta_1\big) + \mathrm{FT}\big(\rho^{-4}\,\tau,\;\rho^{2}\big)\\
+ \mathrm{R}\big(\rho^{-3}\,\eta_1\,\eta_2^{-1}\big) - \mathrm{R}\big(\rho^{-1}\,\tau\,\eta_1\big)\\
+ \mathrm{R}\big(\tau^{-1}\,\eta_2^{-1}\big) - \mathrm{R}\big(\eta_1^{-1}\big)\\
- \mathrm{R}\big(\eta_2\,5^{-1}\big) + \mathrm{R}\big(\tau\,\eta_1^{-1}\,\eta_2\big)\\
+ 11\,\mathrm{R}\big(\tau\big) - \mathrm{R}\big(\tau\,\eta_1\big)\\
= 13[z_0]+[z_1]-10[C_0] .
\end{multline*}
Exceptional pairs $x\mapsto 1-x$: $\rho^{-4}\,\tau\mapsto \rho^{-3}\,\tau^{4}$, $\rho^{-3}\,\eta_1\,\eta_2^{-1}\mapsto \rho^{-2}\,\tau^{2}\,\eta_2^{-1}$, $\rho^{-2}\,\tau\mapsto \rho^{-1}\,\tau^{2}\,\eta_2$, $\rho^{-2}\,\tau^{2}\,\eta_2^{-1}\mapsto \rho^{-3}\,\eta_1\,\eta_2^{-1}$, $\rho^{-1}\,\eta_2^{-1}\mapsto \rho^{-2}\,\eta_1\,\eta_2^{-1}$, $\rho^{-1}\,\tau\mapsto \eta_1^{-1}\,\eta_2^{-1}\,5$, $\rho^{-1}\,\tau\,\eta_1\mapsto \tau^{2}\,\eta_2$, $\rho^{-1}\,\tau^{2}\mapsto \rho^{-2}\,\tau\,\eta_1$, $\tau^{-1}\,\eta_2^{-1}\mapsto \eta_2^{-1}$, $\eta_1^{-1}\,\eta_2^{-1}\mapsto \rho^{2}\,\tau^{-2}\,\eta_1^{-1}\,\eta_2^{-1}$, $\eta_1^{-1}\mapsto \rho\,\eta_1^{-1}$, $\eta_2^{-2}\mapsto \rho^{-2}\,\tau^{-1}\,\eta_1^{2}\,\eta_2^{-2}$, $\eta_2^{-1}\mapsto \tau^{-1}\,\eta_2^{-1}$, $\eta_2\,5^{-1}\mapsto \tau^{2}\,\eta_2\,5^{-1}$, $\tau\,\eta_1^{-1}\,\eta_2\mapsto \rho\,\eta_1^{-2}\,\eta_2^{-1}\,5$, $\tau\,\eta_2^{-1}\mapsto \rho^{2}\,\tau^{-2}\,\eta_1^{-1}\,\eta_2^{-2}\,5$, $\tau\mapsto \rho$, $\tau\,\eta_1\mapsto \rho^{2}$, $\tau\,\eta_1\,\eta_2\,5^{-1}\mapsto \tau^{2}\,\eta_1\,\eta_2^{2}\,5^{-1}$, $\tau^{2}\,\eta_1^{-1}\mapsto \rho\,\tau^{-1}\,\eta_1^{-2}\,\eta_2^{-1}\,5$, $\tau^{2}\mapsto \rho\,\tau\,\eta_2$, $\tau^{2}\,\eta_2\mapsto \rho^{-1}\,\tau\,\eta_1$, $\rho\,\eta_1^{-1}\mapsto \eta_1^{-1}$, $\rho\,\eta_2^{-1}\mapsto \tau^{-1}\,\eta_2^{-2}\,5$, $\rho\mapsto \tau$, $\rho\,\eta_2\,5^{-1}\mapsto \tau^{2}\,\eta_2^{2}\,5^{-1}$, $\rho\,\tau\,\eta_1^{-1}\mapsto \rho^{2}\,\eta_1^{-2}\,5$, $\rho^{2}\,\eta_1^{-2}\mapsto \rho^{2}\,\eta_1^{-2}\,\eta_2$, $\rho^{2}\mapsto \tau\,\eta_1$.

\Needspace{10\baselineskip}
\subsection*{Identity \#8 ($m=31$, $N=19$)}
Generators: $\rho$ the root of $f_2$ in $(0,1)$, $\tau=1-\rho$, $\lambda=\rho^{2}-\rho+1$ (norm $19$), $\eta_1=1+\rho$, $\eta_2=(1+\tau)/\tau$, and the rational prime $7$. Expanding through the pair dictionary and
cancelling yields the exact equality
\begin{multline*}
-5\,\mathrm{FT}\big(\tau,\;\tau\big) + 4\,\mathrm{FT}\big(\eta_1^{-1},\;\tau\,\eta_1\big)\\
- \mathrm{FT}\big(\lambda^{-1}\,\eta_1^{-1},\;\lambda\big) + \mathrm{FT}\big(\lambda,\;\tau\,\lambda^{-1}\big)\\
+ 2\,\mathrm{FT}\big(\rho^{-1}\,\eta_1\,\eta_2^{-1},\;\tau\,\eta_1^{-1}\,\eta_2\big) + 4\,\mathrm{FT}\big(\rho^{-1}\,\tau,\;\rho^{-1}\,\tau\big)\\
+ \mathrm{FT}\big(\tau^{-1}\,\lambda^{-1}\,\eta_2^{-1},\;\lambda\big) + 2\,\mathrm{FT}\big(\rho^{-2}\,\eta_1\,7^{-1},\;\rho^{2}\big)\\
+ \mathrm{FT}\big(\rho^{-1}\,\lambda\,\eta_1^{-1},\;\tau\,\lambda^{-1}\,\eta_1\big) - 2\,\mathrm{FT}\big(\eta_2\,7^{-1},\;\tau^{2}\big)\\
+ \mathrm{FT}\big(\rho^{-1}\,\tau^{-1}\,\lambda\,\eta_2^{-1},\;\rho^{-1}\,\tau\big) - \mathrm{FT}\big(\rho^{-1}\,\eta_1\,\eta_2^{-1},\;\rho^{-1}\,\lambda\,\eta_1^{-1}\big)\\
- \mathrm{FT}\big(\eta_2^{-1},\;\tau^{2}\,\lambda^{-1}\,\eta_2\big) + 4\,\mathrm{FT}\big(\rho^{-2}\,\tau\,\eta_1,\;\eta_2\,7^{-1}\big)\\
+ 4\,\mathrm{FT}\big(\rho^{-1}\,\tau^{2}\,\eta_2,\;\eta_1\,7^{-1}\big) - 2\,\mathrm{FT}\big(\eta_1^{-1}\,\eta_2^{-1},\;\tau\big)\\
+ \mathrm{FT}\big(\rho^{-1}\,\tau,\;\rho\,\lambda^{-1}\big) + \mathrm{FT}\big(\lambda^{-1}\,\eta_1^{-1},\;\tau\,\lambda\,\eta_2^{-1}\big)\\
- 2\,\mathrm{FT}\big(\eta_1^{-1},\;\eta_2^{-1}\big) + 2\,\mathrm{FT}\big(\rho^{-2}\,\tau\,\eta_1\,\eta_2\,7^{-1},\;\tau^{-1}\,\eta_2^{-1}\big)\\
+ \mathrm{FT}\big(\rho^{-2}\,\tau^{2},\;\rho^{-1}\,\eta_1\,\eta_2^{-1}\big) - 2\,\mathrm{FT}\big(\rho^{-1}\,\tau\,\eta_1^{-1}\,\eta_2^{-1}\,7,\;\eta_2\,7^{-1}\big)\\
+ 2\,\mathrm{FT}\big(\rho^{-1}\,\tau^{2}\,\eta_1\,\eta_2\,7^{-1},\;\rho\,\eta_1^{-1}\big) + \mathrm{FT}\big(\rho^{-3}\,\lambda\,\eta_1^{-1},\;\rho^{3}\big)\\
- \mathrm{FT}\big(\tau\,\lambda^{-1}\,\eta_1,\;\tau\,\eta_1^{-1}\,\eta_2^{-1}\big) + \mathrm{FT}\big(\rho^{-6}\,\tau,\;\rho^{3}\,\tau^{-1}\,\lambda\,\eta_1^{-1}\big)\\
+ \mathrm{R}\big(\rho^{-3}\,\tau^{2}\,\eta_1\,\eta_2^{-1}\big) + \mathrm{R}\big(\rho^{-3}\,\tau^{3}\big)\\
- 2\,\mathrm{R}\big(\rho^{-2}\,\eta_1^{2}\,7^{-1}\big) - 4\,\mathrm{R}\big(\rho^{-2}\,\tau\,\eta_1\big)\\
+ \mathrm{R}\big(\rho^{-2}\,\tau^{2}\big) - \mathrm{R}\big(\rho^{-1}\,\tau^{-1}\,\lambda\,\eta_2^{-1}\big)\\
- \mathrm{R}\big(\tau^{-1}\,\lambda^{-1}\,\eta_2^{-1}\big) - \mathrm{R}\big(\tau^{-1}\,\eta_2^{-1}\big)\\
+ 3\,\mathrm{R}\big(\eta_1^{-1}\big) + 30\,\mathrm{R}\big(\tau\big)\\
- 3\,\mathrm{R}\big(\tau\,\eta_1\big) + \mathrm{R}\big(\tau^{2}\,\lambda^{-1}\,\eta_2\big)\\
- \mathrm{R}\big(\rho^{2}\,\tau^{2}\,\eta_1^{-2}\,\eta_2^{-1}\,7\big)\\
= 31[z_0]+[z_1]-24[C_0] .
\end{multline*}
Exceptional pairs $x\mapsto 1-x$: $\rho^{-6}\,\tau\mapsto \rho^{-5}\,\tau^{6}$, $\rho^{-3}\,\lambda\,\eta_1^{-1}\mapsto \rho^{-2}\,\tau^{4}\,\eta_1^{-1}\,\eta_2$, $\rho^{-3}\,\tau^{2}\,\eta_1\,\eta_2^{-1}\mapsto \rho^{-2}\,\tau^{-1}\,\eta_1^{-1}\,\eta_2^{-2}\,7$, $\rho^{-3}\,\tau^{3}\mapsto \rho^{-2}\,\lambda\,\eta_1^{-1}\,\eta_2^{-1}\,7$, $\rho^{-2}\,\eta_1\,7^{-1}\mapsto \rho^{-1}\,\eta_1^{2}\,7^{-1}$, $\rho^{-2}\,\eta_1^{2}\,7^{-1}\mapsto \eta_1\,7^{-1}$, $\rho^{-2}\,\lambda\,\eta_2^{-1}\mapsto \rho^{-1}\,\tau^{-1}\,\eta_1^{-1}\,\eta_2^{-2}\,7$, $\rho^{-2}\,\tau\,\eta_1\mapsto \rho^{-1}\,\tau^{2}\,\eta_2$, $\rho^{-2}\,\tau\,\eta_1\,\eta_2\,7^{-1}\mapsto \rho^{-1}\,\tau^{2}\,\eta_1\,\eta_2^{2}\,7^{-1}$, $\rho^{-2}\,\tau^{2}\mapsto \rho^{-1}\,\eta_1^{-1}\,\eta_2^{-1}\,7$, $\rho^{-1}\,\tau^{-1}\,\lambda\,\eta_2^{-1}\mapsto \eta_1^{-1}\,\eta_2^{-2}\,7$, $\rho^{-1}\,\eta_1\,\eta_2^{-1}\mapsto \tau^{-1}\,\eta_1^{-1}\,\eta_2^{-2}\,7$, $\rho^{-1}\,\lambda\,\eta_1^{-1}\mapsto \eta_1^{-2}\,\eta_2^{-1}\,7$, $\rho^{-1}\,\tau\,\eta_1^{-1}\,\eta_2^{-1}\,7\mapsto \rho^{-2}\,\tau^{2}\,\eta_1\,\eta_2$, $\rho^{-1}\,\tau\,\eta_1^{-1}\mapsto \tau^{2}\,\eta_1^{-1}\,\eta_2^{2}$, $\rho^{-1}\,\tau\mapsto \eta_1^{-1}\,\eta_2^{-1}\,7$, $\rho^{-1}\,\tau^{2}\,\eta_2\mapsto \rho^{-2}\,\tau\,\eta_1$, $\rho^{-1}\,\tau^{2}\,\eta_1\,\eta_2\,7^{-1}\mapsto \rho^{-2}\,\tau\,\eta_1^{2}\,\eta_2\,7^{-1}$, $\tau^{-1}\,\lambda^{-1}\,\eta_2^{-1}\mapsto \tau^{2}\,\lambda^{-1}\,\eta_2^{-1}$, $\tau^{-1}\,\eta_2^{-1}\mapsto \eta_2^{-1}$, $\lambda^{-1}\,\eta_1^{-1}\mapsto \rho^{3}\,\lambda^{-1}\,\eta_1^{-1}$, $\eta_1^{-1}\,\eta_2^{-1}\mapsto \rho^{2}\,\tau^{-2}\,\eta_1^{-2}\,\eta_2^{-2}\,7$, $\eta_1^{-1}\mapsto \rho\,\eta_1^{-1}$, $\eta_2^{-1}\mapsto \tau^{-1}\,\eta_2^{-1}$, $\eta_2\,7^{-1}\mapsto \tau^{2}\,\eta_2^{2}\,7^{-1}$, $\eta_1\,7^{-1}\mapsto \rho^{-2}\,\eta_1^{2}\,7^{-1}$, $\lambda\,\eta_1^{-1}\mapsto \rho\,\tau\,\eta_1^{-1}\,\eta_2$, $\lambda\mapsto \rho\,\tau$, $\tau\,\lambda^{-1}\mapsto \rho^{2}\,\lambda^{-1}$, $\tau\,\lambda^{-1}\,\eta_1\mapsto \rho^{2}\,\lambda^{-1}\,\eta_1^{-1}\,\eta_2^{-1}\,7$, $\tau\,\eta_1^{-1}\,\eta_2^{-1}\mapsto \rho^{3}\,\tau^{-2}\,\eta_1^{-1}\,\eta_2^{-1}$, $\tau\,\eta_1^{-1}\,\eta_2\mapsto \rho\,\eta_1^{-2}\,\eta_2^{-1}\,7$, $\tau\mapsto \rho$, $\tau\,\eta_1\mapsto \rho^{2}$, $\tau\,\lambda\,\eta_2^{-1}\mapsto \rho^{4}\,\tau^{-2}\,\eta_1^{-1}\,\eta_2^{-2}\,7$, $\tau^{2}\,\lambda^{-1}\,\eta_2^{-1}\mapsto \tau^{-1}\,\lambda^{-1}\,\eta_2^{-1}$, $\tau^{2}\,\lambda^{-1}\mapsto \rho\,\lambda^{-1}$, $\tau^{2}\,\lambda^{-1}\,\eta_2\mapsto \rho\,\lambda^{-1}\,\eta_1^{-1}\,\eta_2^{-1}\,7$, $\tau^{2}\mapsto \rho\,\tau\,\eta_2$, $\tau^{2}\,\eta_2\,7^{-1}\mapsto \tau\,\eta_2^{2}\,7^{-1}$, $\rho\,\lambda^{-1}\mapsto \tau^{2}\,\lambda^{-1}$, $\rho\,\eta_1^{-1}\mapsto \eta_1^{-1}$, $\rho^{2}\mapsto \tau\,\eta_1$, $\rho^{2}\,\tau^{2}\,\eta_1^{-2}\,\eta_2^{-1}\,7\mapsto \rho^{3}\,\tau^{-1}\,\lambda\,\eta_1^{-1}$, $\rho^{3}\,\tau^{-1}\,\lambda\,\eta_1^{-1}\mapsto \rho^{2}\,\tau^{2}\,\eta_1^{-2}\,\eta_2^{-1}\,7$, $\rho^{3}\mapsto \rho^{2}\,\tau\,\eta_1^{-1}\,7$.

\Needspace{10\baselineskip}
\subsection*{Identity \#10 ($m=97$, $N=7$)}
Generators as in \#1. Expanding through the pair dictionary and
cancelling yields the exact equality
\begin{multline*}
3\,\mathrm{FT}\big(\lambda,\;\tau\,\lambda^{-1}\big) + 14\,\mathrm{FT}\big(\rho,\;\rho\big)\\
+ 2\,\mathrm{FT}\big(\rho^{-1}\,\tau\,\lambda,\;\rho\,\tau^{-1}\big) - 16\,\mathrm{FT}\big(\tau,\;\rho^{2}\,\tau^{-1}\big)\\
- 14\,\mathrm{FT}\big(\tau,\;\tau\big) + 8\,\mathrm{FT}\big(\tau^{2},\;\rho\big)\\
- 2\,\mathrm{FT}\big(\rho^{-1}\,\tau^{2}\,\lambda,\;\rho\,\lambda^{-1}\big) - 2\,\mathrm{FT}\big(\rho^{-1}\,\tau\,\lambda,\;\tau\,\lambda^{-1}\big)\\
- 4\,\mathrm{FT}\big(\rho^{2}\,\tau^{-1},\;\rho^{2}\big) - \mathrm{FT}\big(\tau,\;\rho^{2}\,\tau^{-3}\,\lambda\big)\\
+ \mathrm{FT}\big(\rho\,\lambda^{-1},\;\rho^{2}\,\tau^{-3}\,\lambda\big) + 2\,\mathrm{FT}\big(\rho^{-1}\,\tau^{2}\,\lambda,\;\rho\,\tau^{2}\big)\\
- \mathrm{FT}\big(\rho^{2}\,\tau\,\lambda^{-1},\;\rho^{4}\,\tau^{-1}\big) + \mathrm{FT}\big(\tau^{4}\,\lambda,\;\rho^{3}\,\lambda^{-2}\big)\\
+ \mathrm{FT}\big(\tau^{4}\,\lambda,\;\rho^{3}\,\tau^{4}\,\lambda^{-1}\big) + \mathrm{R}\big(\rho^{-8}\,\tau^{11}\big)\\
- 2\,\mathrm{R}\big(\rho^{-1}\,\tau^{2}\big) + 48\,\mathrm{R}\big(\tau\big)\\
+ \mathrm{R}\big(\tau^{2}\,\lambda^{-1}\big) - 4\,\mathrm{R}\big(\rho\,\tau^{-1}\big)\\
+ \mathrm{R}\big(\rho^{2}\,\tau^{-2}\,\lambda\big)\\
= 97[z_0]+[z_1]-45[C_0] .
\end{multline*}
Exceptional pairs $x\mapsto 1-x$: $\rho^{-8}\,\tau^{11}\mapsto \rho^{3}\,\tau^{8}$, $\rho^{-1}\,\tau\,\lambda\mapsto \rho^{2}\,\tau^{2}$, $\rho^{-1}\,\tau^{2}\mapsto \tau^{2}$, $\rho^{-1}\,\tau^{2}\,\lambda\mapsto \rho^{-2}\,\tau^{4}$, $\lambda\mapsto \rho\,\tau$, $\tau\,\lambda^{-1}\mapsto \rho^{2}\,\lambda^{-1}$, $\tau\mapsto \rho$, $\tau^{2}\,\lambda^{-1}\mapsto \rho\,\lambda^{-1}$, $\tau^{2}\mapsto \rho^{-1}\,\tau^{2}$, $\tau^{4}\,\lambda\mapsto \rho^{-5}\,\tau^{7}$, $\rho\,\tau^{-1}\mapsto \rho^{2}$, $\rho\,\lambda^{-1}\mapsto \tau^{2}\,\lambda^{-1}$, $\rho\mapsto \tau$, $\rho\,\tau^{2}\mapsto \rho^{-2}\,\tau^{3}$, $\rho^{2}\,\tau^{-3}\,\lambda\mapsto \rho^{4}\,\tau^{-2}$, $\rho^{2}\,\tau^{-2}\,\lambda\mapsto \rho^{3}\,\tau^{-3}$, $\rho^{2}\,\tau^{-1}\mapsto \rho^{2}\,\tau^{-2}$, $\rho^{2}\mapsto \rho\,\tau^{-1}$, $\rho^{2}\,\tau\,\lambda^{-1}\mapsto \rho^{2}\,\tau^{-2}\,\lambda^{-1}$, $\rho^{3}\,\tau^{-3}\mapsto \rho^{2}\,\tau^{-2}\,\lambda$, $\rho^{3}\,\lambda^{-2}\mapsto \rho^{-2}\,\tau^{4}\,\lambda^{-2}$, $\rho^{3}\,\tau^{4}\,\lambda^{-1}\mapsto \rho^{-4}\,\tau^{6}\,\lambda^{-1}$, $\rho^{3}\,\tau^{8}\mapsto \rho^{-8}\,\tau^{11}$, $\rho^{4}\,\tau^{-1}\mapsto \rho^{3}\,\tau^{-4}$, $\rho^{6}\,\lambda^{-1}\mapsto \rho^{4}\,\tau^{-5}\,\lambda^{-1}$.

\endgroup


\section*{Author contributions}
Cetin Hakimoglu-Brown: Conceptualization; Methodology; Formal analysis;
Investigation; Writing---review and editing.

\section*{Funding and competing interests}
No funding was received for this research or the preparation of this
article. The author reports no competing interests.

\Needspace{14\baselineskip}
\section*{Declaration of generative AI use}
ChatGPT (OpenAI) was used for computational assistance, including symbolic
algebra, exact finite searches and interval calculations, and for help
with corrections, checking arguments and drafting revisions to the
mathematical exposition.
Anthropic Claude Opus 5 was used during preparation of an earlier draft
for symbolic and numerical verification, drafting and checking arguments.
Anthropic Claude Opus 5.5 was used in preparing this version, including
drafting Section~\ref{subsec:index} and Appendix~\ref{subsec:closing},
and repeating the finite calculations.
The mathematical evidence is supplied by the derivations, exact
identities and finite calculations reported in this article.
Responsibility for the mathematical statements, interpretation and
references rests with the author.

\end{document}